\documentclass[11pt]{article}
\usepackage{amsmath, amsthm, amssymb, eucal, xypic}

\newcommand{\partentry}[1]{\addtocontents{toc}
{\small\bfseries#1\hfill\thepage\par}}

\makeatletter
\def\@part[#1]#2{%
    \ifnum \c@secnumdepth >\m@ne
      \refstepcounter{part}
      \partentry{\protect\makebox[2em][l]{\thepart}#1}
\else
      \partentry{#1}
    \fi
    {\noindent\normalfont\Large\bfseries\thepart\hspace{1em}#2\par}
    \nobreak
    \vskip 3ex
    \@afterheading}
\def\@spart#1{%
    {\noindent\normalfont\Large\bfseries #1\par} 
     \nobreak
     \vskip 3ex
     \@afterheading}

\renewcommand\section{\@startsection{section}{1}{\z@}%
 				{-3.5ex \@plus -1ex \@minus -.2ex}
				{2ex \@plus.2ex}									
				{\large\bfseries}}
\renewcommand\subsection{\setcounter{subsection}{\value{equation}}
						\stepcounter{equation}
					\@startsection{subsection}{2}{\z@}
                          {1.75ex \@plus.5ex \@minus.2ex}%
                           {-.4em}
                           {\textit}}
\def\@seccntformat#1{\@ifundefined{#1@cntformat}%
	{\csname the#1\endcsname\quad} 
	{\csname #1@cntformat\endcsname}} 
\def\section@cntformat{\thesection.~} 
\def\subsection@cntformat{(\thesubsection)\ }
\renewcommand*\l@section{\mdseries\small\@dottedtocline{1}{1.5em}{2em}}
\makeatother

\numberwithin{equation}{section}
\newtheorem{theorem}[equation]{Theorem}

\newtheorem{lemma}[equation]{Lemma}
\newtheorem{proposition}[equation]{Proposition}
\theoremstyle{definition}
\newtheorem{definition}[equation]{Definition}
\theoremstyle{remark}
\newtheorem{remark}[equation]{Remark}

\newcommand{\cA}{\mathcal{A}}
\newcommand{\cC}{\mathcal{C}}
\newcommand{\cK}{\mathcal{K}}
\newcommand{\calL}{\mathcal{L}}
\newcommand{\cO}{\mathcal{O}}
\newcommand{\calR}{\mathcal{R}}
\newcommand{\frg}{\mathfrak{g}}
\newcommand{\frl}{\mathfrak{l}}
\newcommand{\frk}{\mathfrak{k}}
\newcommand{\frt}{\mathfrak{t}}

\newcommand{\bF}{\mathbb{F}}
\newcommand{\bP}{\mathbb{P}}
\newcommand{\bQ}{\mathbb{Q}}
\newcommand{\bR}{\mathbb{R}}
\newcommand{\bC}{\mathbb{C}}
\newcommand{\bT}{\mathbb{T}}
\newcommand{\bZ}{\mathbb{Z}}
\newcommand{\vep}{\varepsilon}
\newcommand{\xt}{{}^\tau\!}
\newcommand{\SU}{\mathrm{SU}}
\newcommand{\SO}{\mathrm{SO}}
\newcommand{\Or}{\mathrm{O}}
\newcommand{\lbeta}{(\!(\beta)\!)}
\newcommand{\hclass}{\iota_{3}}
\newcommand{\bur}{BU_{r}}
\newcommand{\burq}{BU_{r}^{\bQ}}
\newcommand{\buq}{BU^{\bQ}}

\begin{document}

\title{\textbf{Twisted equivariant $K$-theory with complex coefficients}}

\author{Daniel S.~Freed \and Michael J.~Hopkins \and Constantin Teleman}
\date{May 15, 2006}
\maketitle

\begin{abstract}
\noindent\begin{quote}
Using a global version of the equivariant Chern 
character, we describe the complexified twisted equivariant $K$-theory 
of a space with a compact Lie group action in terms of fixed-point data. 
We apply this to the case of a compact group acting on itself by conjugation, 
and relate the result to the Verlinde algebra and to the Kac numerator at 
$q=1$. Verlinde's formula is also discussed in this context.
\end{quote}
\end{abstract}

\section*{Introduction}
Let $X$ be a locally compact topological space acted upon by a compact 
Lie group $G$. The equivariant $K$-theory $K^*_G(X)$ was defined by Atiyah
and Segal; some foundational papers are \cite{seg1} and \cite{atseg1}. 
\textit{Twisted} versions of $K$-theory, both equivariant and not, have 
recently attracted some attention. The equivariant twistings we consider 
are classified, up to isomorphism, by an equivariant cohomology class 
$[\vep,\tau]\in H_G^1(X;\bZ/2)\times H_G^3(X;\bZ)$. (A \textit{twisting} 
is a representative cocycle for such a class, in some model of equivariant 
cohomology.) For torsion, non-equivariant twistings, the relevant $K$-theory 
was first introduced in \cite{donkar}; subsequent treatments (\cite{ros} 
and, more recently, \cite{ozgerb,at,atseg3}) remove the torsion 
assumption. We recall for convenience the topologists' definition in 
the case $\vep=0$. Because the projective unitary group $\bP\mathrm{U}$ 
has classifying space $K(\bZ;3)$, a class $[\tau]\in H^3(X;\bZ)$ defines 
a principal $\bP\mathrm{U}$-bundle over $X$ up to isomorphism. Let $\bF_X$ 
be the associated bundle of Fredholm endomorphisms. The negative $\xt K(X)$ 
groups are the homotopy groups of the space of sections of $\bF_X$; the 
others are determined by Bott periodicity. In the presence of a group 
action, the equivariant $\xt K$-groups arise from the space of invariant 
sections.\footnote{A technical variation is needed in the equivariant 
case, pertaining to the topology on Fredholm operators \cite{seg2,atseg3}.} 
In this paper, we shall implicitly assume the basic topological properties 
of twisted $K$-theory; their justification is found in the papers mentioned
earlier. However, since $H^1$ twistings get less coverage in the 
literature, we discuss them in more detail in \S\ref{gradingsect}. 

One of the basic results \cite{seg1} of the equivariant theory expresses, 
in terms of fixed-point data, the localisation of $K_G^*(X)$ at prime 
ideals in the representation ring $R_G$ of $G$. The situation simplifies
considerably after complexification, when the maximal ideals in $R_G$ 
correspond to (complex semi-simple) conjugacy classes. Recall that, in the 
non-equivariant case, the Chern character maps complex $K$-theory 
isomorphically onto complex cohomology. The localisation results can 
be assembled into a description of complex equivariant $K$-theory by a 
\textit{globalised Chern character} \cite{slom,atseg2,rosu} supported 
over the entire group. Part I of our paper generalises these results to 
the twisted case: the main result, Theorem \ref{24}, describes 
$\xt K_G^*(X)$, via the \textit{twisted Chern character} (\S\ref{3}), 
in terms of (twisted) equivariant cohomology of fixed-point sets, 
with coefficients in certain equivariantly flat complex line bundles. 
For orbifolds, this is Vafa's \textit{discrete torsion} \cite
{vaf,vafwit}. We give a more detailed outline of Part I in \S1 below.

In Part II, we apply our main result to the $G$-space $X=G$ with the 
conjugation action. For simplicity, we focus on connected groups with 
free $\pi_1$; according to Borel \cite{borel}, these are characterised 
by the property that all centralisers of group elements are connected. 
We do, nonetheless, discuss the $\SO$ groups (Example \ref{46} and 
\S\ref{so3}), to illustrate the phenomena entailed by torsion in $\pi_1$. 

For \textit{regular} twistings, the $R_G$-module $\xt K_G^*(G)$ is 
supported at finitely many conjugacy classes in $G$ (\S\ref{4}). 
Furthermore, $\xt K_G^*(G)$ is an algebra under the \textit{Pontryagin 
product}, and in \S\ref{6} we see how the conjugacy-class decomposition 
diagonalises this product. For \textit{any} compact group $G$, the 
\textit{integral} $\xt K_G(G)$ is isomorphic to the \textit{Verlinde ring} 
of the theory of loop groups, at a certain level related to the twisting 
\cite{turkfreed,freedicm,fht}. For connected $G$, we describe in \S\ref{5} 
an isomorphism between the complexified objects in terms of the Kac 
character formula. The Verlinde algebra also encodes the dimensions of 
spaces of ``non-abelian $\vartheta$-functions"; in \S\ref{7}, we 
incorporate this into the twisted $K$-theory framework. The final 
section illustrates the effect of $H^1$-twistings with a 
detailed study of $G=\SO(3)$. 

\subsection*{Notation:} Given the focus of our paper, \textit{complex 
coefficients will be understood from now on} in cohomology and $K$-theory, 
unless others are indicated. 
\vskip1cm

\part{The twisted equivariant Chern character}

In the first part of a paper, we explain how the equivariant Chern 
character is modified by a twisting, leading to the appearance of
flat line bundles in the coefficients of cohomology.

\section{The idea} 
Let $G$ be a compact Lie groups acting on $X$, and $\tau$ a twisting 
for the $G$-equivariant $K$-theory of $X$. The twisted $K$-theory $\xt 
K_G^*(X;\bZ)$ is a module over the untwisted $K_G^*(X;\bZ)$, and in 
particular over the ring $K_G^0(*;\bZ)=R_G$ of virtual complex 
representations of $G$. Similarly, the complexified $\xt K_G^*(X):=\xt 
K_G^*(X;\bZ)\otimes_\bZ\bC$ is a module over $R_G \otimes_\bZ\bC$, 
henceforth denoted $\bC R_G$. The character identifies $\bC R_G$ with the 
ring of complex-valued polynomial class functions on $G$; the spectrum 
of this ring is the quotient variety $Q:=G_\bC/\!/G_\bC$ of \textit
{geometric invariant theory}, whose points turn out to correspond to the 
semi-simple conjugacy classes in the complexification $G_\bC$, and also 
to the $G$-orbits of \textit{normal}\footnote{Those commuting with their 
hermitian adjoints} elements in $G_\bC$. When $G$ is connected, $Q$ is 
also the quotient $T_\bC/W$ of a complexified maximal torus by the Weyl 
group. 

Our main theorem \eqref{24} identifies, via the Chern character, the 
formal completion $\xt K^*(X)_{q}^\wedge$ at a point $q\in Q$ with 
the \textit{twisted equivariant cohomology} $\xt H_{Z(g)}^*\left(X^g; 
\xt\calL(g)\right)$ of \S\ref{3}. Here, $g$ is a normal group element 
associated to $q$, $Z(g)\subset G$ the unitary part of its centraliser 
and $X^g$ the fixed-point set in $X$ of (the unitary part of) the 
algebraic subgroup generated by $g$. The coefficients live in an 
equivariantly flat line bundle $\xt\calL(g)$ (defined in \ref{211}), which 
varies continuously with $g$. The only novelty is the twisting; in its 
absence, the $\xt \calL(g)$ are trivial and the result is well-known. 
(For torsion twistings on orbifolds, a closely related result was 
independently obtained in \cite{lup}.) Example \ref{25} describes 
these objects in detail for the group $G=\SU(2)$ acting on itself by conjugation.

Because $Q$ is an affine variety, the $\xt K_G^*(X)$ are spaces of global 
sections of sheaves of $\cO$-modules ${}^\tau\cK^*(X)$ over $Q$, obtained 
by Zariski localisation  (see e.g. \cite[II]{hart}). The sheaves 
are coherent if $\xt K_G^*(X)$ is a finite $\bC R_G$-module; for instance, 
this is the case if $X$ is a finite, $G$-equivariant $CW$-complex. 
Passage from $\xt K_G^*(X)$ to sections of ${}^\tau\cK^*$ is a globalisation 
of the \textit{twisted Chern character} $\xt ch:\xt K_G^*(X)\to \xt H_G^*(X)$, 
which we recall in \S\ref{3}. Just as its untwisted version, $\xt ch$ 
only sees the completion of $K$-theory at the augmentation ideal, as a
consequence of the Atiyah-Segal completion theorem \cite{atseg1}. The 
idea of ``repairing" this problem by defining a global Chern character 
over $G$, while implicit in \cite[Prop.~4.1]{seg1} and perhaps folklore, 
was proposed in \cite{bryl} (and carried out for Abelian groups). Finite 
groups had been treated in \cite{slom} (in a different language); they were 
rediscovered in \cite{atseg2}, and a plethora of variations for compact 
groups followed \cite{ezra,duflo,groj,rosu}. When the twisting is null, 
our twisted theorem specialises to an algebraic version of those results, 
expressing the completed stalks $\cK^*(X)_{q}^\wedge$ in terms 
of twisted equivariant cohomology.

Reduction to twisted cohomology is helpful because we can reduce its 
computation to untwisted cohomology via a spectral sequence (\S\ref{3}). 
This is the $\xt ch$-image of the (twisted) Atiyah-Hirzebruch sequence, 
which converges to the augmentation completion of twisted $K$-theory. 
The latter is our completed stalk ${}^\tau \cK^*(X)_1^\wedge$ at the 
identity; it sees only part of the answer,\footnote{Over $\bZ$, the 
augmentation completion sees a bit more, but usually not everything.} 
and can vanish in interesting cases, such as $\xt K_G(G)$ for simple 
groups $G$. In this respect, equivariant $\xt K$-theory for compact $G$ 
can behave like the (untwisted) equivariant theory for a finite group, 
and the ${}^\tau\cK^*$-sheaves can be skyscrapers on finitely many 
conjugacy classes. It is therefore essential to understand the completions 
at all points of the group, and this is what we do.

A $K$-class is determined by its completion at all points. For skyscraper 
sheaves, the resulting description $\xt K_G^*$ is completely satisfactory. 
In general, our twisted Chern character can likely be modified, as in 
\cite{rosu}, to give localisations of the sheaves $\cK^*$ and patching 
maps. (This, however, does not help with computations, because the Atiyah-%
Hirzebruch sequence converges to the completed stalks anyway.) A good 
global target space for $\xt ch$ would seem to be a twisted version of 
equivariant cyclic homology \cite{ezra}, but we will not pursue that here.

\section{The non-equivariant case}\label{3}

In $K$-theory, one first defines the Chern character for a line bundle 
$L$ with first Chern class $x$ as
\[
ch(L)=e^{x},
\]
and then extends this definition to all of $K$-theory using the splitting 
principle. For twisted $K$-theory, the situation is complicated by the fact 
that the description of twisted cohomology classes, which are the values of 
the twisted Chern character $\xt ch$, involves explicit cochains.

In this section, we review the construction of twisted cohomology and of 
$\xt ch$. Readers familiar with Proposition \ref{35}, or willing to accept 
it, may proceed to the next section. We will use the Quillen-Sullivan 
equivalence \cite{quil,sul} between rational homotopy theories of spaces, 
on one side, and commutative differential graded algebras on the other; 
constructions that avoid rational homotopy theory are found in \cite
{ozgerb,freedicm} and the upcoming \cite[II]{atseg3}. 

\subsection{Twisted Cohomology.}\label{sec:twisted-cohomology}
Let $X$ be a space and $\tau\in Z^{3}(X;\bZ)$ a $3$-cocycle (twisting) 
with class $[\tau]\in H^{3}(X;\bZ)$. Let $\beta$ be an indeterminate of 
degree $-2$; for a non-negatively graded vector space $R$, write 
$R\lbeta$ for the graded vector spaces of formal Laurent series:
\[
R\lbeta^{n}=\prod \beta^{i}R^{n+2i}.
\]
We would like to define the twisted cohomology groups ${\xt H}^*(X;\bQ)$ 
to be the cohomology groups of the complex $C^*(X;\bQ)\lbeta$ with 
respect to the differential $(\delta+\beta\tau)$.  The trouble is that 
$(d+\beta\tau)^{2}=\beta^{2}\tau\cup \tau$ which, though cohomologous
to zero, need not actually be zero. We must work instead with a strictly 
commutative differential graded algebra (DGA) ``model" $A^*(X)$ for the 
rational cochains. (This means that $A(X)$ comes equipped with a weak 
equivalence to $C^*(X;\bQ)$, as an $E_{\infty}$ algebra.) If $\tau$ is 
a $3$-cocycle in $A^*(X)$, we define ${\xt H}^*(X)$ to be the cohomology 
of $A^*(X)\lbeta$ with respect to the differential $(d+\beta\tau)$.

For any $\omega \in A^2$, multiplication by $e^{-\beta\omega}$ identifies 
the cohomologies of $d+\beta\tau$ and $d+\beta (\tau +d\omega)$. This shows 
that, up to non-canonical automorphism, $\xt H^*(X)$ depends on the class 
$[\tau]$ alone; the automorphisms come from the multiplicative action of 
$\exp[\beta H^2(X)]$. More generally (as follows from the next proposition), 
any quasi-isomorphism $\varphi: A^*(X) \to B^*$ of differential graded 
algebras gives an isomorphism of $\xt H^*(X)$ with the twisted cohomology 
of $B^*$ with differential $d+\beta\varphi(\tau)$. In this sense, twisted 
cohomology is independent of the DGA model used for cochains. For 
example, when $X$ is a smooth (para-compact, finite dimensional) manifold, 
the real cochains $C^*(X;\bR)$ are quasi-isomorphic to the de Rham complex 
of $X$. In that case, $\xt H^*(X;\bR)$ is canonically isomorphic to the 
cohomology of the de Rham complex $\Omega^*(X)\lbeta$ with respect to the 
differential $(d+\beta\tau)$, where now $\tau$ is a closed $3$-form 
representing the image of $[\tau]$ in $H^{3}(X;\bR)$. Observe that 
$\xt H^*(X;\bR) = \xt H^*(X;\bQ)\otimes\bR$ when $X$ has finite topology.

\begin{proposition}\label{32} 
Filtration by $A$-degree leads to a spectral sequence with $E_2^{p,q}= 
H^p(X)$ for even $q$, vanishing for odd $q$, with $\delta_2=0$ and 
$\delta_3 = \text{multiplication by }[\tau]$, converging weakly to 
$\xt H^{p+q}(X)$. It converges strongly if $H^*(X)$ is finite-dimensional 
in each degree.
\end{proposition}

\begin{proof} Recall \cite[III]{mccleary} that a spectral sequence 
induced from a descending filtration on a complex \textit{converges weakly} 
if $E_\infty=\mathrm{gr}H^*$ and $H^*$ is complete, for the induced 
filtration. It suffices for this that the filtration on the complex 
should be complete and Hausdorff. Strong convergence means, in addition, 
that $H^*$ is Hausdorff; Mittag-Leffler conditions on the $E_r^{p,q}$ 
ensure that. In our case, we let $F^pA^*\lbeta = \prod_{i\ge p} 
A^i\lbeta$; the convergence conditions are clearly satisfied. 
Further, $E_1^{p,2q} = \beta^{-q}A^p$, $E_1^{p,2q+1} = 0$, and, clearly, 
$\delta_1 =d$, $\delta_3 = \beta[\tau]$. 
\end{proof} 

\begin{remark}\label{33}\begin{trivlist}
\item(i) Additional twisting by a flat line bundle on $X$ is allowed; 
the same line bundle will appear in the cohomology coefficients, in 
\eqref{32}.
\item(ii) When a compact group $G$ acts on $X$ and $[\tau] \in H_G^3(X)$, 
we define $^\tau H_G^*(X)$ to be the twisted cohomology of the Borel
construction (homotopy quotient) $EG\times_GX$.
\end{trivlist}
\end{remark}

\subsection{Examples.}\label{34}
(i) Assume that  $\tau$  is a free generator of $H^*(X)$: that is, 
the latter is isomorphic to $R[\tau]$, for some graded ring $R$. 
Then, $\xt H^*(X)=0$. Indeed, $E_4=0$ in (\ref{32}). 
\vspace{-.5ex}\begin{trivlist}
\item(ii) For an example of (i), take $X$ to be a compact connected Lie 
group $G$, and a $[\tau]$ which is non-trivial on $\pi_3(G)$. Even more 
relevant to us, with the same assumptions, is the homotopy quotient $G/G$
for the adjoint action: the equivariant cohomology $H_G^*(G)$ breaks 
up as $H^*(BG)\otimes H^*(G)$, and the twisted cohomology vanishes
again. 
\item(iii) Take $X=T$, a torus, $G=T$ acting trivially on $X$, and let 
$[\tau]\in H^2(BT)\otimes H^1(T)\subset H_T^3(T)$ be defined by a 
non-degenerate bilinear form on the Lie algebra  $\frt$ of $T$. Then, 
$\xt H_T^*(T) = \bC$, in degrees $\dim T \pmod{2}$. More precisely, 
$\xt H^*(T)\cong H^*(T)\lbeta$ (canonically, if $\tau$ is a product 
of cocycles on the two factors), and restriction $\xt H_T^*(T)\to H^*(T)
\lbeta$ lands in top $*$-degree. The $E_3$ term in \eqref{32} is the 
Koszul complex $\Lambda(\frt^*)\otimes \mathrm{Sym}(\frt^*)\lbeta$ 
with differential $\beta\tau$. 
\item(iv) Higher differentials in the spectral sequence involve 
Massey products with $[\tau]$; see \cite[II]{atseg3}. 
\end{trivlist}

\subsection{The twisted Chern character.}\label{sec:twist-chern-char}
Let $\bF$ be a classifying space for $K$-theory, write $\bur$ for the 
component corresponding to (virtual) bundles of virtual rank $r$ 
and $\burq$ for its rationalisation. The Chern character gives an 
identification
\[
\burq\approx\prod_{n=1}^{\infty} K(\bQ,2n).
\]
(We abandon $\beta$ from the notation, for the topological calculations.)
Identify now $BU(1)=K(\bZ,2)$ and choose a model of $\bF$ on which there 
is an action
\[
K(\bZ,2)\times\bF\to \bF
\]
corresponding to tensoring a vector bundle with a line bundle. The
associated bundle
\[
EK(\bZ,2)\times_{K(\bZ,2)}\bF \to BK(\bZ,2)=K(\bZ,3)
\]
is the classifying space for twisted $K$-theory in the following sense: 
view a $3$-cocycle on $X$ as a map $\tau:X\to K(\bZ,3)$.  Then, $\xt 
K(X;\bZ)$ is the group of vertical homotopy classes of lifts in 
\begin{equation}
\xymatrix{
& EK(\bZ,2)\times_{K(\bZ,2)}\bF  \ar[d] \\
	X \ar@{-->}[ur] \ar[r] & K(\bZ,3).
}
\end{equation}
The space $EK(\bZ,2)\times_{K(\bZ,2)}\bF $ breaks up into a disjoint union
\[
EK(\bZ,2)\times_{K(\bZ,2)}\bF = \coprod_{r} 
EK(\bZ,2)\times_{K(\bZ,2)}\bur,
\]
and the group of vertical homotopy classes of lifts in 
\begin{equation}\label{eq:1}
\xymatrix{
& EK(\bZ,2)\times_{K(\bZ,2)}\bur \ar[d] \\
X \ar@{-->}[ur] \ar[r] & K(\bZ,3).
}
\end{equation}
is the summand of $\xt K^{0}(X;\bZ)$ consisting of ``twisted vector bundles'' 
of virtual rank $r$.   We will define the twisted Chern character separately 
on each such summand.

By naturality of $\bQ$-localisation, the group $K(\bZ,2)$ acts on $\burq$, 
and the map $BU\to\burq$ is equivariant for this action.  The twisted 
Chern character is given by composing a section of~\eqref{eq:1} with 
\[
EK(\bZ,2)\times_{K(\bZ,2)}\bur \to 
EK(\bZ,2)\times_{K(\bZ,2)}\burq
\]
to get a lift in
\begin{equation}\label{eq:2}
\xymatrix{
& EK(\bZ,2)\times_{K(\bZ,2)}\burq \ar[d] \\
X \ar@{-->}[ur] \ar[r] & K(\bZ,3).
}
\end{equation}
The real work is then to identify the group of vertical homotopy
classes of lifts in~\eqref{eq:2} with twisted cohomology. For this 
it will be more convenient to work entirely within the context of
$\bQ$-local spaces. Since $\bQ$-localisation commutes with finite
products, the action of $K(\bZ,2)$ on $\burq$ extends canonically to 
an action of $K(\bQ,2)$. There is then a homotopy Cartesian square
\[
\xymatrix{
EK(\bZ,2)\times_{K(\bZ,2)}\burq \ar[r]\ar[d] &
EK(\bQ,2)\times_{K(\bQ,2)}\burq \ar[d]\\
K(\bZ,3)\ar[r]& K(\bQ,3),}
\]
and we can identify the group of vertical homotopy classes of lifts in~\eqref{eq:2}
with the group of vertical homotopy classes of lifts in
\begin{equation}\label{eq:3}
\xymatrix{
& & EK(\bQ,2)\times_{K(\bQ,2)}\burq  \ar[d] \\
X \ar@{-->}[urr] \ar[r] & K(\bZ,3) \ar[r] & K(\bQ,3).
}
\end{equation}

\subsection{Rational homotopy calculation.}
We now use Sullivan's model \cite{sul} for the rational homotopy category 
of spaces in terms of commutative differential graded algebras.\footnote
{Strictly speaking, the DGA theory works well only for \textit{CW} complexes 
of finite type, but in our application the restriction on $X$ can be removed 
\textit{a posteriori}, by a limiting argument.} We replace the diagram of 
spaces~\eqref{eq:3} with a diagram of DGA models, and compute the vertical 
homotopy classes of lifts in~\eqref{eq:3} as homotopy classes of 
factorisations
\begin{equation}
\xymatrix{
& A\left( EK(\bQ,2)\times_{K(\bQ,2)}\burq \right)   \ar@{-->}[dl]\\
A(X) & A\left(K(\bQ,3) \right). \ar[u]\ar[l]
}\label{eq:4}
\end{equation}
In order for the homotopy classes of DGA factorisations to calculate 
the  homotopy classes of lifts, the top algebra must be \textit{free 
nilpotent} over the base $A\left(K(\bQ,3)\right)$; this is the DGA-world 
condition that the vertical arrow should be a cofibration. (Recall that 
$A^*$ is free nilpotent over $B^*$ if it can be constructed from $B^*$ 
``by adding generators one at a time, with no relations". More precisely, 
$A^*$ admits a space $V$ of free, homogeneous algebra generators over $B^*$, 
which is endowed with an increasing filtration $V_0\subset V_1\subset \ldots$ 
satisfying $dV_{k+1} \subset B^*[V_k]$.) 

We now describe the algebras in~\eqref{eq:4}. If $H^*(Y;\bQ)$ is a free 
graded commutative algebra, then a choice of cocycle representatives for 
a basis of $H^*(Y;\bQ)$ realises $H^*(Y;\bQ)$, with differential $d=0$, 
as a DGA model for $Y$. We can therefore take
\begin{align*}
A(K(\bQ,3))&=\bQ[\hclass], \qquad\quad\deg\hclass =3 \\
A(K(\bQ,2))&=\bQ[\omega], \qquad\quad\:\deg \omega =2 \\
A(\burq)&=\bQ[ch_{0}, ch_{1}, ch_{2},\dots]/(ch_{0}=r),
\end{align*}
where the $ch_{n}$ correspond to cocycles of the Chern character $ch = 
ch_{0}+ch_{1}+\dots$. We will henceforth drop the reminder $ch_{0}=r$ 
from the notation. 

We now find a model for the top algebra in \eqref{eq:4}. Let $S^3$ and 
$D^3$ be the $3$-sphere and disk, respectively, and identify $K(\bZ;3)$ 
rationally with $S^3$. We construct $EK(\bQ,2)\times_{K(\bQ,2)}\burq$ 
by gluing the two copies of $D^3\times\burq$ along the common boundary 
$S^2\times\burq$. A DGA model for $S^2$ is $A(S^2)=\bQ[\omega]/
(\omega^2)$; we use $\bQ$ as a model for one copy of $D^3$ and $A(D^3)= 
\bQ[\omega,\hclass]/(\omega^2,\hclass\omega)$ with $d\omega = \hclass$ 
for the other. The restrictions induced by the boundary maps $S^2\subset 
D^3$ become the inclusion of $\bQ$ and the obvious map $\omega\mapsto\omega, 
\hclass\mapsto 0$, respectively.

By multiplicativity of the Chern character, the map
\[
K(\bQ,2)\times \buq\to \buq
\]
is represented by the map of DGA's
\begin{align*}
\bQ[ch_{0},ch_{1},\dots] &\to \bQ[\omega]\otimes \bQ[ch_{0},ch_{1},\dots]
\\
ch &\mapsto e^{\omega}\otimes ch.
\end{align*}
The two boundaries are thus glued using the map $S^2\times\burq \to 
S^2\times \burq$ preserving $\omega$ and sending $ch$ to $e^\omega\cdot 
ch$. As a consequence of the Mayer-Vietoris theorem, a DGA model for the 
glued space $EK(\bQ,2)\times_{K(\bQ,2)}\burq$ is the fibred sum 
\[
A(D^3)[ch_{0}, ch_{1}, ch_{2},\dots]\oplus_{A(S^2)[ch_{0}, ch_{1}, 
	ch_{2},\dots]}\bQ[ch_{0}, ch_{1}, ch_{2},\dots]
\]
with the first map to the base induced by the boundary restriction 
$A(D^3)\to A(S^2)$ and $ch\mapsto ch$, and the second defined by 
$\bQ\to A(S^2)$, $ch\mapsto e^\omega ch$. The fibred sum becomes
isomorphic to
\[
\bQ[\hclass, ch_{0}, ch_{1}, ch_{2},\dots], \quad 
	d(ch) = \hclass ch,
\]
by sending $ch$ above to $e^\omega ch\in A(D^3)[ch_{0}, ch_{1}, ch_{2},
\dots]$, and the diagram~\eqref{eq:4} becomes
\begin{equation}\label{eq:5}
\xymatrix{
& \bQ[\hclass,ch_{0},ch_{1},ch_{2}\dots]   \ar@{-->}[dl]\\
A(X) & \bQ[\hclass] \ar[u]\ar[l].
}
\end{equation}
The bottom arrow corresponds to a $3$-cocycle $\tau\in A^{3}(X)$, and
a factorisation to an element 
\[
\alpha=\alpha_{0}+\alpha_{1}+\dots\in
A^{\text{ev}}(X),\qquad \deg\alpha_{i}=2i
\]
satisfying
\[
(d-\tau)\alpha = 0.
\]
Two factorisations $(\tau,\alpha^{0})$ and $(\tau,\alpha^{1})$ are
homotopic if there exists a homotopy
\[
h(s)=\xi(s)+\eta(s)\,ds \in A(X)\otimes\bQ[s, ds],\qquad \deg s=0
\]
satisfying
\begin{align*}
(d-\tau)h &=0, \\
h(0) =\alpha^{0},\quad &h(1) =\alpha^{1}.
\end{align*}
These conditions equivalent to the existence of $\zeta\in
A^{\text{odd}}(X)$, satisfying
\[
(d-\tau)\zeta = \alpha^{1}-\alpha^{0}.	
\]
Indeed, given $\xi+\eta\,ds$ we take $\zeta = -\int_{0}^{1}\eta(s)\,ds$, 
and given $\zeta$, we define
\[
\xi+\eta\,ds = s\,(d-\tau)\zeta - \zeta\, ds.
\]
It follows that the group of factorisations is
\[
\ker(d-\tau)/\mathrm{Im}(d-\tau) = {\xt H}^*(X;\bQ).
\]

\subsection{Summary.} We have defined, for spaces $X$ of finite type, 
a Chern character map 
\[
ch: \xt K^{0}(X;\bZ) \to {\xt H}^{\text{ev}}(X;\bQ).
\]
inducing an isomorphism $\xt K^{0}(X;\bQ) \cong {\xt H}^{\text{ev}}(X;\bQ)$.
We leave to the reader to check that this agrees with the usual Chern 
character when $\tau=0$, and that it is multiplicative in the sense 
that the following diagram commutes
\[
\xymatrix{
\xt K^{0}(X;\bZ)\times {}^{\tau'}\!K^{0}(X;\bZ) \ar[r]\ar[d] & 
{}^{\tau+\tau'}\!H^{0}(X;\bQ)\ar[d] \\
{\xt H}^{\text{ev}}(X;\bQ)\times {}^{\tau'}\!H^{\text{ev}}(X;\bQ) \ar[r] 
& {}^{\tau+\tau'}\!H^{\text{ev}}(X;\bQ).
}\]
For reference, we record this in the following Proposition. 

\begin{proposition}\label{35} There exists a functorial twisted Chern
character $\xt ch\!:\xt K^*(X;\bQ)\to \xt H^*(X;\bQ)$, which is a module 
isomorphism over the Chern isomorphism $ch\!:K^*(X;\bQ)\to H^*(X;\bQ)$.
\end{proposition}

\noindent The finite type assumption on $X$, a standard requirement in 
the DGA theory, is removed by taking inverse limits over the finite 
sub-complexes, when $X$ is an arbitrary \textit{CW} complex. Note, 
however, that the isomorphism $K^*(X;\bQ)\cong \xt K^*(X;\bZ) 
\otimes \bQ$ requires $X$ to be a finite complex.

\section{The main result}\label{2}
As a motivating example, we review the case of finite groups. Recall 
that cohomology and $K$-theory have complex coefficients.

\begin{theorem}[\cite{atseg2}, Theorem 2]\label{21} 
For finite $G$, we have 
\[
K_G^*\left(X\right)\cong \bigoplus\nolimits_q 
	K^*\left(X^{g(q)}\right){}^{Z(g(q))}, 
\]
the sum ranging over conjugacy classes $q$ with representatives $g(q)$.
\end{theorem}

\noindent The isomorphism extends the identification of $\bC R_G$ with 
the ring of class functions by the character. For each $g\in G$, the 
restriction to $X^g$ of a $G$-vector bundle $V$ on $X$ splits into 
eigen-bundles $V^{(\alpha)}$ under the fibre-wise $g$-action, and \eqref
{21} sends $V$ to the complex linear combination $\sum{\alpha\cdot V^{
(\alpha)}}$.  More abstractly, let $\langle g \rangle \subset G$ be the 
subgroup generated by $g$; a class in $K_G^*(X)$ restricts to a 
$Z(g)$-invariant one in $K_{\langle g \rangle}^*(X^g) \cong R_{\langle 
g \rangle}\otimes K^*(X^g)$, and taking the trace of $g$ on the first 
factor gives the desired class in $K^*(X^g)$. The Chern character gives 
an isomorphism		
\begin{equation}\label{23}
  ch:K_G^{0,1}(X) \xrightarrow{\sim} \bigoplus\nolimits_q 
	H_{Z(g(q))}^{ev,odd}\left(X^{g(q)}\right).
\end{equation}
We will generalise this to include twistings and arbitrary compact Lie 
groups.

\subsection{Notation.} Let $G$ be any compact Lie group. Since we take an 
algebraic route, we must discuss the complex group $G_\bC$, even though 
all information will be contained in the unitary part $G$. Let $g\in G_\bC$ 
be a normal element, generating an algebraic subgroup $\langle g \rangle_
\bC \subset G_\bC$, with centraliser $Z_\bC$. Normality of $g$ 
ensures that these subgroups are the complexifications of their intersections 
$\langle g\rangle$ and $Z=Z(g)$ with $G$, that $\langle g \rangle $ is 
topologically cyclic, and that $Z$ is its commutant in $G$. These are 
consequences of the fact that normal elements are precisely those contained 
in the complexification of some \textit{quasi-torus} in $G$ (a subgroup 
meeting every component of $G$ in a translate of a maximal torus). 

For a $G$-space $X$, call $X^g\subset X$ the $\langle g \rangle$-fixed point 
set and $\xt K_G^*(X)_g^\wedge$ the formal completion of the complex twisted 
$K$-groups at the conjugacy class of $g$. (These are the completed stalks 
of our $\cK$-sheaves at $g$.)

\subsection{The flat line bundles.} The twisted analogue of \eqref{23} 
involves some local systems $\xt\calL(g)$ that we now describe. (However, 
see also \S\ref{216} below for their heuristic interpretation.) Let $Y$ 
(soon to be $X^g$) be a $Z$-space on which $\langle g \rangle$ acts trivially. 
As $\langle g \rangle$ is topologically cyclic, $H^3(B\langle g \rangle;
\bZ)=0$, so $[\tau] \in H_Z^3(Y;\bZ)$ defines, in the Leray sequence 
\[
H_{Z/\langle g\rangle}^p \left(Y; H^q(B\langle g\rangle;\bZ)\right)
\Rightarrow H_Z^{p+q}(Y;\bZ), 
\]
a class 
\begin{equation}\label{210}			
\mathrm{gr}[\tau]\in H_{Z/\langle g \rangle}^1
	\left(Y; H^2(B\langle g\rangle;\bZ)\right) = 
	\mathrm{Hom}\left(H_1^{Z/\langle g \rangle}(Y;\bZ)\times 
	\langle g \rangle;\bT \right).
\end{equation}

\begin{definition}\label{211} 
The holonomy of the $Z/\langle g \rangle$-equivariant flat complex 
line bundle $\xt\calL(g)$ over $Y$ is the value of the (complexified) 
homomorphism $\mathrm{gr}[\tau]$ at $g$.
\end{definition}

\noindent This defines $\xt\calL(g)$ up to non-canonical isomorphism. 
We can also construct it `on the nose' in two equivalent ways, as 
follows.

\subsection{Construction.}\label{212} 
Choose a twisting $\tau$ representing $[\tau]$ in the form of a projective 
Hilbert bundle $\bP_Y \to Y$, with projective-linear lifting of the $Z
$-action. Over every point of $Y$, a central extension of the group 
$\langle g \rangle $ by the circle is defined from the action of 
$\langle g \rangle \subset Z$ on the fibres of $\bP_Y$. The complexified 
extensions define a principal $\bC^\times$-bundle over $Y\times\langle 
g \rangle_\bC$, and we call $\xt\calL$ the associated complex line bundle.
\begin{trivlist}\itemsep0ex
\item(i) $\xt\calL(g)$ is the restriction of $\xt\calL$ to $Y\times\{g\}$. 
Flatness is seen as follows. Since $\langle g\rangle$ is topologically 
cyclic, the extensions are trivial, but not canonically split.  The 
splittings, however, form a discrete set, any two of them differing by a 
homomorphism $\langle g\rangle_\bC\to \bC^\times$. Following a splitting 
around a closed loop defines such a homomorphism, and its value at $g$ 
gives the holonomy of $\xt\calL (g)$. 
\item(ii) The $1$-dimensional characters of the centrally extended
$\langle g \rangle$'s which restrict to the natural characters of the 
centres form a principal bundle $p\!:\tilde Y\to Y$ with fibre $H^2(B
\langle g \rangle;\bZ)$, equivariant under $Z/\langle g \rangle$. The 
isomorphism class of $\tilde Y$ is $\mathrm{gr}[\tau]$ in (\ref{210}). 
Evaluation at $g$ defines a homomorphism $H^2(B\langle g \rangle;\bZ)\to 
\bC^\times$; thereunder, the line bundle $\xt\calL(g)$ is associated to 
$\tilde Y$. 
\end{trivlist}

\subsection{Twisted Chern isomorphism for compact $G$.} 
We now generalise \eqref{23} to the case of a finite \textit{CW}-complex 
$X$ a  with $G$-action. The result extends to other $G$-spaces, but may 
require some adjustment, depending on which version of $K$-theory is used. 
(For instance, with proper supports, $\xt K_G^*(X;\bZ)\otimes_{R_G}
(\bC R_G)_g^\wedge$ replaces the completion.) 

\begin{theorem}\label{24} 
The formal completions $\xt K_G^*(X)_g^\wedge$ at $g$ are isomorphic to 
the twisted cohomologies $\xt H_{Z(g)}^*\left(X^g;\xt \calL(g)\right)$.
\end{theorem}

\noindent The reader may now wish to consult Example \ref{25} at the 
end of this section. 

\begin{proof}
The theorem assembles the two propositions that follow. 
\end{proof}

\begin{proposition}\label{26} 
The restriction $\xt K_G^*(X)_g^\wedge \to \xt K_Z^*(X^g)_g^\wedge$ 
is an isomorphism.
\end{proposition}

\begin{proposition}\label{215} 
There is a natural isomorphism $\xt K_Z^*(X^g)_g^\wedge \cong 
\xt H_Z^*\left(X^g;\xt\calL(g)\right)$.  
\end{proposition}

\noindent The isomorphism will be constructed in the proof. Its indirect 
nature is largely responsible for the length of the argument.

\begin{remark}\label{27} 
\begin{trivlist}\itemsep0ex
\item (i) Proposition \ref{26} resembles \cite[Prop.~4.1]{seg1}, but is 
not quite equivalent to it. As our proof shows, \eqref{26} still holds if 
\textit{\'etale localisation} (Henselisation) is used instead of completion, 
but fails for Zariski localisation, even when $X$ is a point: the fraction 
field of $T$ differs from that of $Q$, which is the Weyl-invariant subfield.
\item(ii) When $g=1$, $\xt\calL(g)$ is the trivial line bundle, and 
Proposition \ref{215} results by combining the completion theorem \cite
{atseg1} with the twisted Chern isomorphism of \S2. In this sense, we 
are generalising the completion theorem to the case of a twisting and  
a central element; the new feature is the appearance of the flat line 
bundle. However, we will use the original result in the proof.
\end{trivlist}\end{remark}

\begin{proof}[Proof of (\ref{26}).]
We start with a point $X=*$, when $\xt K_G^0(*;\bZ)$ is the module, 
also denoted $\xt R_G$, of \textit{$\tau$-projective representations}. 
More precisely, the class $[\tau] \in H_G^3(*;\bZ)$ defines an isomorphism 
class of central extensions of $G$ by the circle group $\bT$. Fixing such 
an extension $\xt G$---which can be viewed as a cocycle representative 
for $[\tau]$---allows one to define the abelian group of those virtual 
representations of $\xt G$ on which the central circle acts naturally. 
This is the topologist's definition\footnote{In the $K$-theory of 
C*-algebras, this is a theorem.} of $\xt K_G^0(*;\bZ)$; it is an 
$R_G$-module, under tensor product. The complexification of $\xt G$ 
defines an algebraic line bundle $\xt\cC h^0$ over $G_\bC$, carrying a 
natural lifting of the conjugation action. Its fibres are the lines 
$\xt\calL(g)$ over $X=*$. Its invariant direct image to $Q$ is a 
torsion-free sheaf ${}^\tau\cK^0$. (This need not be a line bundle, 
because the centralisers may act non-trivially on some fibres). 
Characters of  $\tau$-projective representations are invariant sections 
of $\xt\cC h^0$, and examining class functions on $\xt G$ identifies 
the complexification $\xt K_G^*(*)$ with the space of invariant 
sections of $\xt\cC h^0$, which is also the space of sections of 
${}^\tau\cK^0$. 

When $\tau =0$, we are asserting that the completed ring $(\bC R_G)_g
^\wedge$ of class functions is isomorphic, under restriction from $G$ 
to $Z$, to $(\bC R_Z)_g^\wedge$. This is true because $Z_\bC$ is a local 
(\'etale) slice at $g$ for the adjoint action of $G_\bC$. For general 
$\tau$, the sheaf $\xt\cC h^0$ on $Z_\bC$ is the restriction of its 
$G_\bC$-counterpart, so the two direct images ${}^\tau\cK^0$ on $Q$, 
coming from $G_\bC$ and from the local slice $Z_\bC$, agree near $g$.

To extend the result to a general $X$, it suffices, by a standard argument,
to settle the case of a homogeneous space $G/H$. In that case, $\xt K_G^*
(G/H) = \xt K_H^*(*)$. The conjugacy class of $g$ in $G_\bC$ meets $H_\bC$ 
in a finite number of classes,\footnote{The morphism from $Q_H$ to $Q_G$ 
is affine and proper, hence finite. (Properness is seen by realising the 
$Q$'s as the spaces of \textit{unitary} conjugacy classes of normal elements.)} 
for which we can choose normal representatives $h_i = k_i^{-1}gk_i$, with 
$k_i\in G$, and get a natural isomorphism
\begin{equation}\label{28}			
\xt K_G^*\left(G/H\right)_g^\wedge \cong \bigoplus\nolimits_i
\xt K_H^*(*)_{h_i}^\wedge \cong \bigoplus\nolimits_i\xt 
K_{k_iHk_i^{-1}}^*(*)_g^\wedge.
\end{equation}
A coset $kH\in G/H$ is invariant under $\langle g\rangle$-translation 
iff $k^{-1}\langle g \rangle k\in H$. This holds precisely when 
$k^{-1}gk\in H_\bC$; thus, $k^{-1}gk=hh_ih^{-1}$, for some $h\in H_\bC$ 
and a unique $i$. As $k^{-1}gk$ and $h_i$ are both normal, we can assume 
$h\in H$. Hence, $khk_i^{-1}\in Z$, and $k\in Z k_iH$, for a unique $i$. 
The fixed-point set of $\langle g \rangle$ on $G/H$ is then the disjoint 
union over $i$ of the subsets $Z\cdot k_iH$. These are isomorphic, as 
$Z$-varieties, to the homogeneous spaces $Z\left/{Z\cap k_iHk_i^{-1}}
\right.$. From here,
\begin{equation}\label{29}				
\xt K_Z^*\left((G/H)^g\right)\cong \bigoplus\nolimits_i \xt 
K_{Z\cap k_iHk_i^{-1}}^*(*),
\end{equation}
and, as $Z\cap k_iHk_i^{-1}$ is the centraliser of $g$ in $k_iHk_i^{-1}$,
equality of right-hand sides in (\ref{28}) and (\ref{29}), after 
completion at $g$, follows from the known case when $X$ is a point.   
\end{proof}

\subsection{The twisting $\tau'$.} \label{213}
For the proof of Proposition \ref{215}, we need one more construction.
Note that the pull-back projective bundle $\bP_{\tilde Y}:= p^*\bP_Y$ to 
$\tilde Y$ in Construction \ref{212} has a \textit{global} $g$-eigenspace 
decomposition: the (globally defined) projective sub-bundles labelled by 
the characters of $\langle g\rangle$ split the (locally defined) Hilbert 
space fibres underlying $\bP_{\tilde Y}$. The twisted $K$-groups defined 
by any \textit{non-empty} sub-bundle are naturally isomorphic to those 
defined by $\bP_{\tilde Y}$. Now, the inclusion of the $g$-trivial 
eigen-bundle $\bP'_{\tilde Y}\subset \bP_{\tilde Y}$ gives a reduction 
of $p^*\tau$ to a $Z/\langle g \rangle$-equivariant \textit{twisting} 
$\tau'$. (We can ensure non-emptiness by tensoring $\bP_Y$ with the 
regular representation of $Z$.) This refines the fact that the class 
$p^*[\tau] \in H_Z^3(\tilde Y;\bZ)$ comes from a distinguished 
$[\tau']\in H_{Z/\langle g \rangle}^3(\tilde Y;\bZ)$.

\begin{lemma}\label{214} The composition ${}^{\tau'}\!K_{Z/\langle g
\rangle}^*(\tilde Y)\to {}^{\tau'}\!K_Z^*(\tilde Y)\buildrel {p_!} \over 
\longrightarrow \xt K_Z^*(Y)$ is an isomorphism. (Here, $K$-theory has 
proper supports along $p$.)
\end{lemma}

\begin{proof} A Mayer-Vietoris arguments reduces us to the case when 
$Y$ is a homogeneous $Z$-space, or equivalently (after replacing $Z$ 
by the isotropy subgroup) to the case when $Y$ is a point. The lemma 
then says $\xt R_Z={}^{\tau'}\! K_{Z/\langle g \rangle}^0(\tilde Y)$. 
In view of our definition of $\tau'$, that simply expresses the 
eigenspace decomposition of $\tau$-twisted $Z$-representations under 
the central element $g$.  
\end{proof}

\begin{remark} \begin{trivlist}\itemsep0ex
\item (i) The result and its proof apply to \textit{integral} 
$K$-theory.
\item (ii) A similar argument shows that the inverse map is the 
$g$-invariant part of $p^*$.
\end{trivlist}\end{remark}
 
\begin{proof}[Proof of (\ref{215})]
We will use the isomorphism in Lemma \ref{214} to construct a map $\xt 
K_Z^*(Y) \to \xt H_Z^*\left(X^g;\xt\calL(g)\right)$, and complete the 
$R_Z$-modules at $g$ in two steps in order to produce the desired isomorphism.
We outline the argument before checking details.

First, we complete at the augmentation ideal of the ring $\bC R_{Z/
\langle g\rangle}$. By the completion theorem \cite{atseg1}, this converts  
${}^{\tau'}\!K_{Z/\langle g\rangle}^*(\tilde Y)$ into the ${}^
{\tau'}\!K$-theory of the homotopy quotient of $\tilde Y$ by $Z/\langle 
g\rangle$. The twisted Chern character equates that with the $\tau'$-twisted, 
$Z/\langle g \rangle$-equivariant cohomology of $\tilde Y$. (In both cases, 
we use proper supports along $p$ over every finite sub-complex of $BZ/
\langle g\rangle$.) The map $p$ re-interprets this cohomology as the twisted 
$Z/\langle g \rangle$-equivariant cohomology of $Y$ with coefficients 
in a flat vector bundle $\cC$. By construction of $p$ \eqref{212}, the 
fibres of $\cC$ are the spaces of polynomial sections of $\xt\calL$ over 
$\langle g\rangle$, and are free $\bC R_{\langle g \rangle}$-modules 
of rank one. We then complete $\cC$ at $g$ and identify the result with 
the bundle $\xt\calL(g)[[\frl]]$ of $\xt\calL(g)$-valued power series on
the Lie algebra $\frl$ of $\langle g \rangle$. Finally, we identify the 
twisted $H^*_{Z/\langle g \rangle}\left(Y;\xt\calL(g)[[\frl]]\right)$ with 
the desired $\xt H_Z^*\left(Y; \xt\calL(g)\right)$ by re-assembling $\frl$ 
and the Lie algebra of $Z/\langle g \rangle$ into that of $Z$. 

We now construct the cohomological counterpart of the maps in \eqref{214} 
in the Cartan model for equivariant cohomology. Splitting the Lie algebra 
of $Z$ into $\frl$ and a Lie complement $\frk$ decomposes the equivariant 
$3$-forms as
\begin{equation}\label{components}
\Omega^3_Z(Y) = \frl^*\otimes\Omega^1_{Z/\langle g\rangle}(Y)
	\oplus\Omega^3_{Z/\langle g\rangle}(Y).
\end{equation} 
For degree reasons, $\Omega^1_{Z/\langle g\rangle}(Y) =\Omega^1(Y)^Z$. 
We use a $3$-form $\tau$ to represent our twisting in the Cartan model. 
(See Remark \ref{220}.ii below for the relation to the projective bundles.) 
This $\tau$ decomposes as $\tau^{(1)} \oplus\tau^{(3)}$. Since 
$[\tau^{(1)}]= \mathrm{gr} [\tau]$ is the defining $1$-class for the 
covering space $\tilde Y$, we have $p^*\tau^{(1)} = d\omega$ for some 
$\frl^*$-valued, $Z$-invariant function $\omega$ on $\tilde Y$.

We use $p^*\tau^{(3)}$ as a de Rham model for the twisting $\tau'$. 
Choose a basis $\zeta_a$ of $\frk$, with dual basis $\zeta^a$. The 
twisted $ch$ then identifies the augmentation completion of ${}^{\tau'}
\!K_{Z/\langle g\rangle}^*(\tilde Y)$ with the cohomology of 
\[
\left\{\Omega^*_{Z/\langle g\rangle}(\tilde Y),\, 
	d + \zeta^a\iota(\zeta_a) + p^*\tau^{(3)}\wedge\right\},
\]
where the mod $2$ grading of the Cartan complex allowed us to drop 
$\beta$ from the coefficients. Projecting to $Y$ leads to
\begin{equation}\label{incomplete}
\left\{\Omega^*_{Z/\langle g\rangle}(Y;\cC), 
	\,\nabla + \zeta^a\iota(\zeta_a) + \tau^{(3)}\wedge\right\},
\end{equation}
with the natural flat connection $\nabla$ on $\cC$. 

We now identify the fibre-wise completion of $\cC$ at $g$ with 
$\xt\calL(g)[[\frl]]$ as follows. View a section $\varphi$ of $\cC^\wedge_g$ 
as a function on the formal neighbourhood of $\tilde Y\times\{g\}$ 
inside $\tilde Y\times\langle g\rangle$, invariant for the action 
of the character group $H^2(B\langle g\rangle;\bZ)$ (\ref{212}.ii). 
We send $\varphi$ to 
\begin{equation}\label{defpsi}
\psi: (y,\xi)\mapsto \exp(-\langle\omega(\tilde y) | \xi\rangle)\cdot
	\varphi\left(\tilde y, g\cdot\exp(\xi)\right)	
\end{equation}
where $\xi\in\frl$ and $\tilde y$ is any lift of $y\in Y$. The defining 
properties $d\omega = p^*\tau^{(1)}$ for $\omega$ and (\ref{212}.ii) for 
$\tilde Y$ ensure that $\psi$ is independent of the lift $\tilde y$. By the 
same formula for $d\omega$, $\nabla\varphi$ maps in this way to $(\nabla_g 
+ \tau^{(1)}\wedge)\psi$, where $\nabla_g$ is the flat connection on 
$\xt\calL(g)$. The complex \eqref{incomplete} therefore becomes, upon 
completion,  
\[
\left\{\Omega^*_{Z/\langle g \rangle}\left(Y;\xt\calL(g)[[\frl]] \right), 
\, \nabla_g + \zeta^a\iota(\zeta_a) + \tau\wedge\right\},
\] 
which our splitting $\mathrm{Lie}(Z) = \frl\oplus\frk$ identifies with 
the Cartan model for $\xt H_Z^*\left(Y;\xt\calL(g)\right)$.
\end{proof}

\begin{remark}\label{220} \begin{trivlist}\itemsep0ex
\item (i) The exponential pre-factor in \eqref{defpsi} is the multiplication 
by $e^{-\beta\omega}$ of \S\ref{sec:twisted-cohomology}, which gives the 
isomorphism between the $p^*\tau$- and $p^*\tau^{(3)}$-twisted, 
$Z$-equivariant cohomologies of $\tilde Y$. This matches the reduction 
to $\bP'_{\tilde Y}$ in Lemma \ref{214}. 
\item (ii) Our $\omega$ can be shifted by an $\frl^*$-valued constant. 
This ambiguity appears because we have not reconciled the projective model 
for twistings of Lemma \ref{214} with the de Rham model used to define $ch$. 
To do that, we must choose compatible $2$-forms $\alpha$ on the projective 
bundles, with the $d\alpha$ equal to the lifts of the $\tau$. The 
$\frl^*$-component of $\alpha$ on $\bP'_{\tilde Y}$ then turns out to be 
lifted from the base $\tilde Y$, giving the choice of $\omega$ compatible 
with our Chern character. 
\end{trivlist}\end{remark}

\subsection{Meaning of $\xt\calL(g)$.} \label{216} 
The following heuristic argument for Proposition \ref{215}, inspired by
the case of finite groups, illuminates the appearance of the flat line 
bundles. Translation by $g$ on conjugacy classes defines an automorphism 
$T_g$ of the algebra $\bC R_Z$; this sends an irreducible representation 
$V$ of $Z$, on which $g$ must act as a scalar $\xi_V(g) \in\bC$, to 
$\xi_V(g) \cdot V$. We can lift $T_g$ to an intertwining automorphism 
of the $R_Z$-module $K_Z^*(Y)$ by decomposing vector bundles into $\langle 
g \rangle$-eigenspaces, as in Theorem~\ref{21}. We would like to assert 
the following twisted analogue:
\vskip1.6ex
\emph{$T_g$ lifts to an intertwining $R_Z$-module isomorphism 
$\xt K_Z^*(Y)\cong\xt K_Z^*\left(Y;\xt\calL(g)\right)$.}
\vskip1.6ex
\noindent This would identify the completion of $\xt K_Z^*(Y)$ at $g$ 
with that of $\xt K_Z^*\left(Y;\xt\calL(g)\right)$ at $1$, leading to the 
desired equivariant cohomology via Chern characters and the completion 
theorem \cite{atseg1}. The map from left to right would send a bundle $V$ 
to the linear combination $\sum {\alpha\cdot V^{(\alpha)}}$, defined from 
its $g$-eigenspace decomposition. However, $\langle g \rangle $ acts 
projectively, so its eigenvalues $\alpha$ are flat sections of $\xt 
\calL(g)$, rather than constants.

However, no construction of equivariant $K$-theory with coefficients 
in flat line bundles having the ``obvious" properties exists; our 
work-around was Lemma \ref{214}. 

\subsection{Example.}\label{25} 
Let $G=\SU(2)$, acting on $X=\SU(2)$ by conjugation. 
Then, $G_\bC=\mathrm{SL}(2;\bC)$, and $Q$ may be identified with the 
affine line, with coordinate $q$, as follows: the conjugacy class of the 
matrix $g=\mathrm{diag}(\lambda,\lambda ^{-1})$ corresponds to the point 
$q=\lambda+\lambda ^{-1}$, for $\lambda \in \bC^\times$. These matrices 
form a complexified maximal torus $T_\bC$, and the Weyl group $S_2$ 
interchanges $\lambda$ and $\lambda^{-1}$. The closed interval $[-2,2]$ 
is the image of $\SU(2)$ in $Q$. The unitary centraliser $Z(g)$ is the 
maximal torus $T$, unless $q=\pm 2$, in which case it is the entire 
group $G$. The algebraic group $\langle g \rangle_\bC\subseteq G_\bC$ 
generated by $g$ is $T_\bC$, and its unitary part $\langle g\rangle$ equals 
$T$, unless $\lambda$ is a root of unity, in which case $\langle g \rangle 
=\langle g \rangle_\bC$ is the finite cyclic group generated by $g$. The 
fixed-point set $X^g=T$, unless $\lambda =\pm 1$, in which case $X^g=G$. 
The twistings are classified by $H_G^3(G;\bZ)\cong\bZ$; we focus on the 
case $[\tau] \ne 0$. If $g\ne \pm\mathrm{I}$, the flat line bundle $\xt\calL
(g)$ over $T$ has holonomy $\lambda^{2[\tau]}$ (see \S\ref{4}); so the 
cohomology in Thm.~\ref{24} is nil, unless $\lambda^{2[\tau]} = 1$. 

At $\lambda_k:=\exp(k\pi i/[\tau])$ ($k=1,\ldots,[\tau]-1$), the line 
bundle is trivial, but the differential which computes the twisted 
cohomologies in (\ref{24}) is not, and is described as follows (\ref
{34}.iii). Write $H_T^*(T)\cong H^*(BT)\otimes H^*(T)\cong \left.\bC[[u,
\theta]] \right/\theta^2$, where $u$ is the coordinate on the Lie algebra 
of $T$ (so $\lambda =\lambda_k\exp(u)$ is a local coordinate on $T$); the 
twisted cohomology group is then identified with the cohomology of the 
differential $2[\tau]\cdot u\theta $, and is the line $\bC\theta$. 

On the other hand, at $g=\pm\mathrm{I}$, $\xt\calL(g)$ is trivial, but now 
$H_G^*(G)$ is identified with the sub-ring of Weyl invariants in $H_T^*(T)$, 
which is $\bC[[u^2,u\theta]]\left/(u\theta)^2\right.$; the twisting 
differential $2[\tau]\cdot u\theta$ has zero cohomology, so there is 
no contribution from those points if $[\tau] \ne 0$. All in all, we find
that $\xt K_G^0(G)=0$, while $\xt K_G^1(G) = \bC^{[\tau]-1}$, supported 
at the points $q_k=\lambda_k+\lambda_k^{-1}$ in $Q$. 

When $[\tau] = 0$, a similar discussion shows that the $\cK^{0,1}$ 
are locally free of rank one over $Q$; for a generalisation of this 
result, see \cite{brylz}.

\section{Gradings, or twistings by $H^1(\bZ/2)$}\label{gradingsect}
We now discuss the changes to Theorem \ref{24} in the presence of an
additional $K$-theory twisting, also called \textit{grading}, $\vep \in 
H_G^1(X;\bZ/2)$; the end result is (\ref{A14}). The ideas and definitions 
should be clear to readers familiar with the $K$-theory of graded $C^*$ 
algebras, as in \cite[\S14]{black} (see Remark \ref{A8} below); but we 
discuss group actions and their Chern character sheaves ${}^\vep\cK^*$ in 
more detail. 

In this section, $G$ is any compact Lie group. Also, save for (\ref
{A5}--\ref{A6r}) and (\ref{A11}--\ref{A14}), which involve characters 
and localisation, our discussion applies to \emph{integral} $K$-theory.

\subsection{If $X$ is a point.} A class in $H_G^1(*;\bZ/2)$ corresponds 
to a \textit{$\bZ/2$-grading} of $G$, a homomorphism $\vep:G\to\{\pm 1\}$. 
The fibres $G^\pm$ of $\vep$ are conjugation-stable unions of components 
of $G$. A \textit{graded representation} is a $\bZ/2$-graded vector space 
with a linear action of $G$, where even elements preserve, and odd ones 
reverse the grading. We use the notation $M^+\ominus M^-$, the superscript 
indicating the eigenvalue of the \textit{degree operator $\deg$}; the 
reason for this (purely symbolic) ``direct difference" notation will be 
clear below. A \textit{super-symmetry} of such a representation is an 
odd automorphism, skew-commuting with $G$; representations which admit 
a super-symmetry are called \textit{super-symmetric}. 

\begin{definition}\label{A1} 
${}^{\vep,\tau}R_G = {}^{\vep,\tau }K_G^0(*;\bZ)$ is the abelian group 
of finite-dimensional, $\tau$-twisted, graded representations, modulo 
super-symmetric ones.
\end{definition}

The sum of any graded representation with a degree-reversed copy of
itself is super-symmetric; because of this, flipping the grading acts as
a sign change on $K^0$, and so (\ref{A1}) defines an abelian group, 
not just a semi-group. \textit{Restriction of graded representations}, 
from $G$ to $G^+$, sends $M^+\ominus M^-$ to the virtual representation 
$M^+-M^-$. Restriction has a right adjoint \textit{graded induction} 
functor from $\xt R_{G^+}$ to ${}^{\vep,\tau}R_G$. Both are $R_G$-module 
maps. When $\vep$ is not trivial, graded induction identifies ${}^{\vep,
\tau} R_G$ with the co-kernel of ordinary restriction from $\xt R_G$ to 
$\xt R_{G^+}$, and graded restriction with the kernel of ordinary induction 
from $\xt R_{G^+}$ to $\xt R_G$. This is further clarified in (\ref{A4}) 
below. The last description shows that $^{\vep,\tau}R_G$ is a free 
abelian group.

We define the graded $K^1$ by an implicit Thom isomorphism. The group
$G\times\bZ/2$ carries a $\{\pm 1\}$--valued $2$-cocycle $\kappa$, lifted 
by $\vep\times\mathrm{Id}$ from the Heisenberg extension of $\bZ/2\times
\bZ/2$. Call $\nu$ the grading on $G\times\bZ/2$ defined by the second 
projection.

\begin{definition}\label{A2} 
${}^{\vep,\tau}R_G^1= {}^{\vep,\tau}K_G^1(*;\bZ):={}^{\vep\nu,\tau+
\kappa}R_{G\times \bZ/2}$. 
\end{definition}

\begin{remark} The Thom isomorphism for graded $C^*$-algebras 
identifies ${}^{\vep,\tau}K_G^1(*)$ with the $\xt K^0$-group of the 
\textit{graded crossed product} of $(G;\vep)$ with the graded rank $1$ 
Clifford algebra $C_1$. This product is the graded convolution algebra 
${}^{\tau+\kappa}C_*(G\times\bZ/2)$, with product grading $\vep\nu$; 
the cocycle $\kappa$ stems from the anti-commutation of odd elements. 
\end{remark}

\begin{proposition}\label{A3}
${}^{\vep,\tau}R_G^1\cong \xt R_G\left/\mathrm{Ind}\xt R_{G^+}\right.$
\end{proposition}

\begin{proof} We claim that an $\vep\nu$-graded, $(\tau+\kappa)$-twisted 
representation $M^+\ominus M^-$ of $G\times\bZ/2$ is determined up to 
canonical isomorphism by $M^+$, which is a $(\tau+\kappa)$-twisted 
representation of the even, $(\mathrm{Id},\vep)$-diagonal copy of $G$ in 
$G\times\bZ/2$. Indeed, denoting by ${}^\vep\bC$ the sign representation 
of $G$, the isomorphism $M^-\cong M^+\otimes{}^\vep\bC$ defined by the 
element $(1,-1)\in G\times\bZ/2$ fixes the other summand and generates 
the full group action. Super-symmetric representations are those for which 
$M^+$ is isomorphic to its $\vep$-twist: these are precisely the ones 
induced from $G^+$. Finally, note that, while $\kappa$ is not trivial 
as a $\bZ/2$-extension over the diagonal $G$ (its class is $\vep^2$), 
the resulting $\bT$-extension is trivialised by a choice of $\sqrt{-1}$. 
This leads to our description of ${}^{\vep,\tau}R_G^1$.
\end{proof}

\subsection{The exact sequences.} Inclusion of $G$ as the even subgroup 
of $G\times \bZ/2$ defines graded restriction and induction maps, $\mathrm
{Res'}: {}^{\vep,\tau}R_G^1\to \xt R_G$, $\mathrm{Ind'}: \xt R_G\to 
{}^{\vep,\tau}R_G^1$. In our concrete description of $R^1$, they send 
$M^+$ to $M^+-M^+ \otimes{}^\vep\bC$ and $M$ to $M^+:=M$. When $\vep\neq 
1$, these maps assemble to the two exact sequences generalising (\ref{A3}),
\begin{equation}\begin{split}\label{A4}		
0\to {}^{\vep,\tau}R_G^1\xrightarrow{\mathrm{Res'}} \xt R_G
\xrightarrow{\mathrm{Res}} \xt R_{G^+}\xrightarrow{\mathrm{Ind}} 
{}^{\vep,\tau}R_G\to 0, \\
0\leftarrow {}^{\vep,\tau}R_G^1\xleftarrow{\mathrm{Ind'}}\xt 
R_G\xleftarrow{\mathrm{Ind}} \xt R_{G^+}\xleftarrow
{\mathrm{Res}} {}^{\vep,\tau}R_G\leftarrow 0.
\end{split}\end{equation}
When $\vep= 1$, $\xt R_{G^+}$ is replaced by $\xt R_G\oplus\xt R_G$; 
$\mathrm{Res}$ is the diagonal and $\mathrm{Ind}$ the difference map. 
The uniform entry is $\xt K_G^0(\{+,-\})$, $G$ acting via $\vep$.

\subsection{Super-characters.}\label{A5}
Recall from \S\ref{2} that $\tau$ defines a $G_\bC$-equivariant 
line bundle $\xt\cC h$ over $G_\bC$, whose invariant sections are spanned 
by the characters of irreducible,  $\tau$-projective representations. 
Anti-invariant sections are those transforming under the character $\vep$ 
of $G_\bC$; the \textit{super-character} of a graded representation, 
$g\mapsto Tr(g\deg)$, is an example. The super-character of a graded 
$G$-representation is supported on $G^+$, because odd group elements 
are off-diagonal in the $M^\pm$-decomposition, while the super-character
of a $(\vep\nu,\tau +\kappa)$-twisted representation of $G\times \bZ/2$ 
lives on $G^-\times \{-1\}$, because, with $S=(1,-1)$, $Tr(g\deg) = 
Tr(S\cdot g\deg \cdot S^{-1}) = -Tr(g\deg)$, if $g\in G^+$. The following 
is clear from the exact sequences (\ref{A4}).

\begin{proposition}\label{A6} The $\bC R_G$-module $^{\vep,\tau}K_G^0(*)$ 
is isomorphic to the space of anti-invariant algebraic sections of 
$\xt\cC h$ on $G_\bC^+$, while $^{\vep,\tau}K_G^1(*)$ is isomorphic 
to the space of invariant sections of $\xt\cC h$ over $G_\bC^-$. Both 
isomorphisms are realised by the super-character. \qed 
\end{proposition}

\begin{remark}\label{A6r}
In terms of the ``odd line" $\bC^-$, $^{\vep,\tau }R_G$ is the graded 
module of skew-invariant sections of $\xt\cC h\otimes\bC^-$ over $G$. 
Over $G^+$, $\bC^-$-valued functions are odd, so skew invariance means 
anti-invariance. Over $G^-$, the sections are even, and skew invariance 
means invariance.
\end{remark}

\subsection{General $X$.} On a space $X$, $\vep\in H_G^1(X;\bZ/2)$ 
defines a real $G$-line bundle $\calR$ over $X$, with unit interval 
bundle $D$. The statement below is really a twisted Thom isomorphism 
theorem; we thank G.~Segal for the suggestion of using it as a 
definitional shortcut.

\begin{definition}\label{A7} The groups ${}^{\vep ,\tau }K_G^*(X;\bZ)$
are the relative $\xt K_G^{*-1}$-groups of $(D,\partial D)$.
\end{definition}

\begin{remark} \label{A8}
The boundary $\partial D$ is a double $G$-cover $p:\tilde X
\to X$, and we can give a $C^*$-friendly description of ${}^{\vep,\tau}
K_G^*(X)$ as the $\xt K^{*+1}$-groups of the crossed product $G\ltimes 
C^*(\tilde X)$, graded by the deck transformation. These $K$-groups do 
\textit{not} usually have a naive description in terms of graded projective 
modules, as in the case of a group algebra (\ref{A1},\ref{A2}); see 
\cite{black} for the correct definitions. 
\end{remark}

From the vanishing of the first and third Stiefel-Whitney classes of 
the bundle $\calR\oplus\calR$ over $X$, the Thom isomorphism allows us to 
identify the doubly $\vep$-twisted $K$-groups ${}^{2\vep,\tau}K_G^*$ 
with $\xt K_G^*$, but a choice of sign is needed (in the Spin structure). 
The various \textit{integral} $K$-groups are related by two six-term exact 
sequences analogous to (\ref{A4}), involving the double cover $\tilde X$. 
For notational convenience, we omit the twisting $\tau$ and the integral
coefficients,  present everywhere. The first one is the six-term sequence 
for $D$ and its boundary; the second follows from the first, by replacing 
$X$ by the pair $(D,\partial D)$ and using the triviality of $p^*\vep$ on 
$\tilde X$.

\begin{equation}\label{A9}	
\begin{matrix}		
{K_G^0(X)}&\xrightarrow{\;p^*}&{K_G^0(\tilde X)}&
\xrightarrow{\;{}^\vep p_!}&{{}^\vep K_G^0(X)}\cr
\uparrow &{}&{}&{}&\downarrow \cr
{{}^\vep K_G^1(X)}&\xleftarrow{{}^\vep p_!\;}
&{K_G^1(\tilde X)}&\xleftarrow{\,p^*\,} &{K_G^1(X)}
\end{matrix}\end{equation}
\vspace{2ex}
\begin{equation}\label{A8i}
\begin{matrix}
{K_G^0(X)}&\xleftarrow{\,p_!\:} &{K_G^0(\tilde X)}&
\xleftarrow{{}^\vep p^*\:} &{{}^\vep K_G^0(X)}\cr
\downarrow &{}&{}&{}&\uparrow \cr
{{}^\vep K_G^1(X)}&\xrightarrow{\;{}^\vep p^*} 
& {K_G^1(\tilde X)}&\xrightarrow{\;p_!} &{K_G^1(X)}
\end{matrix}
\end{equation}
Moreover, denoting by $\alpha$ the deck transformation of $\tilde X$, we have 
\begin{align}
 p_!p^* &= \left(1+\calR\right)\otimes_\bR,
&p^*p_! &= \left(1+\alpha^*\right),\nonumber\\
{}^\vep p_!{}^\vep p^* &= \left(1+\calR\right)\otimes_\bR, 
&{}^\vep p^*{}^\vep p_! &= (1-\alpha ^*).  
\end{align}
Noting that the operations $\alpha^*$ and $\calR\otimes$ both square to $1$, 
we can decompose, modulo $2$-torsion, all the $K$-groups in 
(\ref{A9}) by eigenvalue, and conclude the following.  

\begin{proposition}\label{A10} 
Modulo $2$-torsion, we have 
${}^{\vep,\tau}K_G^*(X)\cong \xt K_G^*(\tilde X)^-\oplus\xt K_G^{1-*}(X)^-$ 
for the integral $K$-groups; the superscripts indicate the $(-1)$-eigenspaces 
of $\alpha^*$, resp.~$\calR\otimes_\bR$, as appropriate. 
Similarly, $K_G^*(X)\cong \xt K_G^*(\tilde X)^+\oplus
{}^{\vep,\tau}K_G^{1-*}(X)^-$. \qed 
\end{proposition}

\subsection{Localisation.} \label{A11}
Finally, we have the following graded version of the results in \S\ref{2}. 
The first proposition is proved just as Proposition \ref{26}, by reduction 
to homogeneous spaces, in view of our description \eqref{A4} of the graded representation rings.

\begin{proposition}\label{A12} The restriction 
${}^{\vep,\tau}K_G^*(X)_g^\wedge \to {}^{\vep,\tau}K_Z^*(X^g)_g^\wedge$ 
is an isomorphism. \qed
\end{proposition}

\subsection{}\label{A13}
An element $g\in G$ is called \textit{even} or \textit{odd}, over a 
component $Y$ of its fixed-point set $X^g$, according to the sign of 
$\vep(g)$. An even $g$ fixes $p^{-1}(Y)$, an odd one switches the fibres. 
Let ${}^{\vep,\tau}\calL(g)$ denote the flat line bundle $\xt\calL(g)\otimes_\bR\calR$. 

\begin{theorem}\label{A14} \begin{trivlist}
\item(i) If $g$ is even over $Y$, ${}^{\vep,\tau}K_Z^*(Y)_g^\wedge 
\cong \xt H_Z^*\left(Y;{}^{\vep,\tau}\calL(g) \right)$. 
\item(ii) If $g$ is odd over $Y$, ${}^{\vep,\tau}K_Z^*(Y)_g^\wedge \cong 
\xt H_Z^{*+1}\left(Y;\xt\calL(g)\right)$. 
\end{trivlist}\end{theorem}

\begin{proof} Regardless of the parity of $g$, we have from \eqref{A10} 
and \eqref{A12}
\begin{equation}\label{A15}
{}^{\vep,\tau}K_Z^*(Y)_g^\wedge \cong \xt K_Z^*\left(p^{-1}(Y)\right)
_g^{\wedge-}\oplus \xt K_Z^{*+1}(Y)_g^{\wedge-}.
\end{equation}
When $g$ is even, Theorem \ref{24} gives natural isomorphisms 
\begin{align*}		
\xt K_Z^*\left(p^{-1}(Y)\right)_g^\wedge &\cong \xt
		H_Z^*\left(p^{-1}(Y);p^*\xt\calL(g) \right), \\
\xt K_Z^{*+1}(Y)_g^\wedge &\cong \xt 
		H_Z^{*+1}\left(Y;\xt\calL(g)\right).
\end{align*}
The $(-1)$-eigenspace for $\alpha^*$ in the first term is $\xt H_Z^*\left
(Y; {}^{\vep,\tau}\calL(g)\right)$. Tensoring with $\calR$ on the second term 
acts as the identity, because $ch\calR = 1$ and $g$ acts trivially on its 
fibres; so the second summand in (\ref{A15}) is nil, in this case.\footnote
{Vanishing also follows from surjectivity of the complexified maps 
$p_!$ on the completions in (\ref{A8i}), which is clear when they are 
identified, by Thm.~\ref{24}, with the $p_*$ maps in cohomology.}

On the other hand, when $g$ is odd, it acts freely on $p^{-1}(Y)$, so
$\xt K_Z^*(p^{-1}(Y))_g^\wedge =0$, by (\ref{26}); while $\calR\otimes_\bR$
acts as $(-1)$ on $\xt K_Z^{*+1}(Y)_g^\wedge$, since $g=-1$ on $\calR$. 
\end{proof}

\part{ A special case: $\xt K_G^*(G)$}

In the remainder of the paper, we study a compact, connected Lie group 
$G$ acting on itself by conjugation.

\section{Reduction to the maximal torus}\label{4}
Let $[\tau] \in H_G^3(G;\bZ)$ be an integral class which restricts 
trivially to $H^3(T)$, the maximal torus. Call $\tau$ \textit{regular} 
when its restriction to $H^1(T)\otimes H^2(BT)\subset H_T^3(T)$ has full 
rank, when viewed as a linear map $H_2(BT)\to H^1(T)$. We shall assume 
this to be so. The last map defines, after tensoring with $\bT$, an 
isogeny $\lambda: T\to T^\vee$ to the dual torus $T^\vee $, whose kernel 
is a finite subgroup $F\subset T$. A class $[\vep]\in H_G^1(G;\bZ/2)$ 
defines, by restriction to $H^1(T;\bZ/2)$, a $2$-torsion point $\vep^\vee
\in T^\vee$, and we let $F_\vep$ denote its $\lambda$-inverse. Interpreting 
points in $T^\vee$ as isomorphism classes of flat line bundles on $T$, 
we observe the following.

\begin{proposition}\label{41} The flat line bundles $\xt\calL(t)$ and ${}^
{\vep,\tau}\calL(g)$ over $T$ (\ref{211}, \ref{A13}) are classified 
by $\lambda(t)$ and $\vep^\vee\cdot\lambda(t)$, respectively. \qed 
\end{proposition}

Using Theorem \ref{24}, this determines $\xt K_G^*(G)$ in some important 
special cases.

\begin{theorem}\label{42}\begin{trivlist}\itemsep0ex
\item(i) For the trivial action of $T$ on itself, ${}^\tau\cK^*(T)$ is a 
skyscraper sheaf with one-dimensional stalks, supported at the points 
of $F$, in dimension $\dim T \pmod{2}$.
\item(ii) Let $\pi_1G$ be free. For the conjugation action of $G$ on 
itself, ${}^\tau\cK^*(G)$ and ${}^{\vep,\tau}\cK^*(G)$ are the skyscraper 
sheaves on $Q=T_\bC/W$ with $1$-dimensional stalks, in dimension 
$\dim G \pmod{2}$, supported on the \emph{regular} Weyl orbits in $F$ and 
$F_\vep$, respectively.
\end{trivlist}\end{theorem}

\begin{proof} This follows from Examples \ref{34}.ii,iii together
with part (i) of the following lemma.   \end{proof} 

\begin{lemma}\label{43} Let $G$ be a connected compact group. 
\vspace{-.8ex}
\begin{trivlist}\itemsep0ex
\item(i) If $\pi_1G$ is free, the centraliser of any element is connected.
\item(ii) For any connected $G$ and $g\in G$, $\pi_0Z(g)$ is naturally a 
subgroup of $\pi_1G$. 
\item(iii) If $g$ is regular, $\pi_0Z(g)$ embeds into the Weyl 
group of $G$.
\end{trivlist}
\end{lemma}

\begin{proof}\begin{trivlist}
\item(i) Let $G'\subset G$ be the commutator subgroup and $G^{ab}$ the 
abelianisation of $G$. The short exact sequence $0\to \pi_1G'\to \pi_1G
\to \pi_1G^{ab}\to 0$ shows that $\pi_1G'$ is the torsion subgroup of 
$\pi_1G$; hence, $G'$ is simply connected. The neutral component $Z_1$ 
of the centre of $G$ surjects onto the quotient $G^{ab}$, because the 
two Lie algebras are isomorphic and $G^{ab}$ is connected; so $G=G'Z_1$. 
Translating by $Z_1$ reduces our assertion to the case of $G'$; but 
this is a result of Borel \cite{borel}.
\item(ii) Write $G={\tilde G}/\pi$, where the central subgroup $\pi$ in 
the covering group $\tilde G$ is isomorphic to the torsion subgroup of 
$\pi_1G$, and $\pi_1\tilde G$ is free. Any $Z(g)$-conjugate of a lifting 
$\tilde g$ of $g$ is another lift of $g$, and their ratio defines a 
homomorphism from $Z(g)\left/{Z(g)_1}\right.$ to $\pi$. If $z\tilde 
gz^{-1}=\tilde g$, then $z$ lifts to an element of the centraliser of 
$\tilde g$ in $\tilde G$; but the latter is connected, by (i), so $z$ 
must lie in $Z(g)_1$, and our homomorphism is injective.
\item(iii) By regularity, $Z(g)_1$ is the unique maximal torus containing 
$g$, and is normalised by $Z(g)$; so $Z(g)\left/{Z(g)_1}\right.$ embeds 
into the Weyl group. 
\end{trivlist}\end{proof}

\begin{remark}\label{44}  
Any connected $G$ contains regular elements $g$ for which $\pi_0 Z(g)$ 
is the full torsion subgroup of $\pi_1G$. When $G$ is simple, such a $g$ 
is $\exp(\rho/c)$; here, $\rho$ is the half-sum of the positive roots, 
$c$ the dual Coxeter number of $\frg$ and $\frt$ is identified with 
$\frt^*$ via the \emph{basic} inner product, which matches long roots 
with short co-roots. For general $G$, one can use the product of these 
points in the simple factors. 
\end{remark}

\subsection{Torsion in $\pi_1$.}\label{45}
Singular conjugacy classes cannot contribute to $\xt K_G(G)$, because the 
relevant twisted cohomology vanishes when $Z_1$ is not a torus (Example 
\ref{34}.ii). The contribution of a regular point $f\in T$ is the invariant 
part, under $Z\!:=Z(f)$, of $\xt H_T^*\left(Z;\xt\calL(f)\right)$. Each 
component $Z_w$ of $Z$ is a translate of $T$, labelled by an element 
$w\in\pi_1G\subset W$ (\ref{43}.iii), and can only contribute to cohomology 
if the flat $T$-line bundle $\xt\calL(f)$ on $Z_w$ is trivial. If $\tau$ 
is regular, this happens at isolated points $f$, and then $\xt H_T^*\left(
Z_w;\xt\calL(f)\right)$ is one-dimensional and confined to top degree. 
All in all, we get one line from $Z_w$ when the flat line bundle $\xt\calL(f)$ 
is $Z$-equivariantly isomorphic to $\det(\frt^w)$, the volume form of the 
quotient $Z_w/T$ (by conjugation); else, we get no contribution. Thus, 
${}^\tau \cK^*(G)$ is still a sum of skyscraper sheaves, but a clean 
general expression for the rank of its stalks seems difficult. To say more, 
let us assume that there exists a point $z\in Z_w$ whose centraliser $Z(z)$ 
meets all the components of $Z$. Construction~\ref{212} defines an action of 
$Z(z)$ on the fibre at $z$ of $\xt\calL(f)$: specifically, this is the 
conjugation action of $Z(z)$ on the fibre over $f$ of the central 
extension defined by the restriction of $[\tau]$ to $H^3_{Z(z)}(z;\bZ)$. 
This action must agree with that of $\pi_0Z$ on $\det(\frt^w)$, if $Z_w$ is 
to contribute to ${}^\tau \cK^*$. 

\subsection{Example.} \label{46} We describe the ingredients in \S\ref{45}
for the $\SO$ groups; the centraliser actions on $\calL$ are justified at the 
end. The  maximal torus $T$ in $\SO(2n)$ and $\SO(2n\!+\!1)$ 
is $\SO(2)^n$, identified with $\bR^n/\bZ^n$ so that the eigenvalues of 
$\boldsymbol\xi$ in the standard representation are $\exp(\pm2\pi i\,\xi_p)$. 
The positive Weyl chamber for $\SO(2n)$ is defined by $|\xi_1|\leq \xi_2\leq 
\ldots\leq \xi_n$, with the extra condition $\xi_1\geq 0$ for $\SO(2n\!+\!1)$. 
The conjugacy classes correspond to the points therein with $\xi_n\leq 1/2$; 
for $\SO(2n)$, we need the extra identification $\xi_1\leftrightarrow-\xi_1$ 
when $\xi_n=1/2$. The equivariant $H^3$ group is\footnote{$\SO(4)$ has an 
extra $\bZ$.} $\bZ/2\oplus\bZ$, split by restriction to $H^3_\SO(e)\cong
\bZ/2$, and the isogeny for the generator of the free part is the standard 
identification of $T$ with its dual torus.\footnote{On Spin groups, we get 
a four-fold cover $\tilde{T} \to\tilde{T}^\vee$: the integer lattice is 
$\langle\pm\mathbf{e}_p \pm\mathbf{e}_q\rangle$.} The points $\boldsymbol\xi$ 
for which $\xt\calL (\boldsymbol\xi)$ is trivial on $T$ are the roots of unity 
of order $k$, the free part of $[\tau]$. 

In $\SO(2n)$, regular points with $\xi_1=0,\:\xi_n=1/2$ have an extra 
component in the stabiliser; conjugation by $w$, the product of the first 
and last coordinate reflections in $\bR^{2n}$, switches the signs of $\xi_1$ 
and $\xi_n$, and acts on the fibre of $\calL$ over $T$ as $\pm(-1)^k$, the 
initial sign matching the torsion part of $[\tau]$. The quotient $Z_w/Z$ 
of the odd component under conjugation is the torus $(\bR/\bZ)^{n-2}$, 
missing the first and last factors of $T$, and $\xt\calL(\boldsymbol\xi)$ 
has holonomy $k\cdot\xi_p$ along $\mathbf{e}_p$. Let $g_p$ be the torus
element corresponding to $\mathbf{e}_p/2$; the centraliser of $w$ in 
$Z$ is $\langle w, g_1,g_n\rangle\cong\{\pm1\}^3$. In the $\tau$-extension, 
$g_n$ will $\pm(-1)^k$-commute with $w$ and $(-1)^k$-commute with $g_1$; 
so $w$ acts on $\xt\calL$ as $\pm(-1)^k$ and $g_1$ as $(-1)^k$.

For $\SO(2n\!+\!1)$, we find a new component in the centraliser of regular 
points with $\xi_n=1/2$; the product $w$ of the last two coordinate reflections
on $\bR^{2n+1}$ switches the sign of $\xi_n$, and acts on $\calL$ over $T$ by 
the same $\pm(-1)^k$. Over the component $Z_w$, with quotient $(\bR/\bZ)^
{n-1}$, the holonomy of $\xt\calL$ along $\mathbf{e}_p$ is again trivial 
precisely over the points with $k\cdot\xi_p\in\bZ$ ($1\leq p\leq n-1$). 
The centraliser of $w$ is $\langle w,g_n\rangle\cong\{\pm1\}^2$, $w$ 
acts on $\calL$ as $\pm(-1)^k$ and $g_n$ acts trivially. 

To summarise, the exceptional $f$-points for $G=\SO(2n)$ have $\xi_1=0$, 
$\xi_n=1/2$, and $\xi_p\in\frac{1}{k}\bZ$ for $1<p<n$; both $T$ and $Z_w$ 
contribute one line to $\cK^n$ if $\tau= (+,\mbox{even})$, and not otherwise. 
For $\SO(2n\!+\!1)$, the exceptional points have $\xi_n=1/2$ and $\xi_p
\in\frac{1}{k}\bZ$ for $1\leq p<n$; $T$ contributes to $\cK^n$ iff $\tau= 
(-,\mbox{even})$, while $Z_w$ contributes to $\cK^{n-1}$ for twistings 
$(+,\mbox{even})$ and $(-,\mbox{odd})$. We hasten to add, though, that a 
meaningful treatment of the orthogonal groups must include the $H^1
$-twistings, as in \S\ref{so3}.

\emph{The twistings.} In the even case, the $\tau$-twisting of the centraliser 
of $g_n$ is determined by restricting to the $\SO(4)$ subgroup acting on the 
first two and last two coordinates of $\bR^{2n}$. Indeed, the four possible 
central extensions are determined by the signs of the commutators of $w$ with
$g_1$ and $g_n$ ($g_1,g_n$ always commute). We claim that $w$ $\pm$-commutes 
with $g_1$ and $\pm(-1)^k$-commutes with $w$. The $\pm$ factor is checked 
directly in the Spin group. To detect the other factor, we must reach a bit
ahead to \S\ref{57}, and embed $Z$ inside the loop group $LG$, identifying 
$T$ with the constant loops and $w$ with $w'=w\cdot\gamma_n$, the straight 
loop from $1$ to $\mathbf{e}_n$. Restricting the $\tau$-central extension 
of $LG$ to this copy of $Z$ captures the restriction of $[\tau]$ to $H^3_{Z}
(g_n;\bZ)$. Now, in the extended $LG$, $\gamma_n$ commutes with $g_1$ and 
$(-1)^k$-commutes with $g_n$, as desired. As $w$ is swapped with $g_n$ by a conjugation in $\SO$ which also swaps $g_1$ with $w\cdot g_1\cdot g_n$, it 
follows that, in the centraliser $\langle w,g_1,g_n\rangle\cong\{\pm1\}^3$ 
of $w$, the commutators are as asserted earlier.

In the odd case, we can detect the twistings in the $\SO(3)$ subgroup on 
the last three coordinates: the centraliser $\Or(2)$ therein is extended by 
the class $\pm(-1)^k$. This can be checked as above, or by a Mayer-Vietoris 
argument as in \S\ref{so3}. 

\subsection{}Instead of a case-by-case listing, we shall indicate in the 
next section a direct relation of the sheaves $\cK$ to representations of 
the loop group of $G$ in terms of the \textit{Kac character formula}, in
the $q\to1$ limit. We shall confine the detailed discussion to the case 
when $\pi_1$ is free; however, we wish to state the general result, and 
this requires a slight variation of Theorem \ref{42}.

Writing $G\cong\tilde{G}/\pi$, as in the proof of (\ref{43}.ii), gives an 
isomorphism $\xt K^*_G(G)=\xt K^*_{G\times\pi}(\tilde{G})$. The new ${}^\tau
\cK^*$-sheaves live on $Q\times\pi$, but they merely unravel the old ones, 
so that the $Z_w$-contribution is supported on $Q\times\{w\}$, and all stalks 
are now one-dimensional. The extra $\pi$ in the stabiliser acts freely 
on $\tilde{Z}_w$, with quotient $Z_w$, so that the $\pi$-trivial holonomy 
condition on the line bundle upstairs matches the trivial holonomy 
condition downstairs. To each $f\in T$ and $w\in\pi_0Z$, the twisting 
$\tau$ assigns the holonomy representation $\pi_1^Z(Z_w)\to\bT$ of $\xt\calL(f)$. 
We reformulate the discussion in \S\ref{45} as follows.

\begin{proposition}\label{47}
The new sheaves ${}^\tau\cK^*(\tilde{G})$ have one-dimensional stalks,
supported at the pairs $f\in Q^{reg}$, $w\in\pi_0Z(f)$ where the holonomy
representation agrees with $\det(\frt^w)$. \hfill \qed
\end{proposition}

\section{The Kac numerator and $\xt K$-class of a loop group 
representation}\label{5}
Our description of $\xt K_G^*(G)$ leads to a concrete if intriguing 
isomorphism with the complexified space of \textit{positive energy 
representations} (PERs) of the loop group $LG$ at a \emph{shifted} 
level, via their distributional characters. The discussion that follows 
is quite crude, as it ignores the \emph{energy} action on representations, 
so that we only see the Kac numerator at $q=1$; detecting the $q$-powers 
requires the rotation-equivariant version of $\xt K_G^*(G)$ \cite{fht}.

\subsection{Level of a central extension.} The representations that concern 
us are projective, and the relevant central extensions of $LG$ by $\bT$ 
turn out to be classified by their topological \textit{level} in $H_G^3
(G; \bZ)$. This arises form the connecting homomorphism $\partial$ in the 
exponential sequence for group cohomology with smooth coefficients,
\begin{equation}\label{51}			
H_{LG}^2(\bR)\to H_{LG}^2(\bT)\xrightarrow{\:\partial\:}H_{LG}^3(\bZ)
\to H_{LG}^3(\bR),
\end{equation}
and the identification, for connected $G$, of $H_{LG}^3(\bZ)$ with
$H_G^3(G;\bZ)$. (For connected $G$, $BLG=LBG$, and the latter is the
homotopy quotient $G/G$ under the adjoint action.) When $G$ is semi-simple, 
the outer terms vanish (\cite{ps}, Ch.~14) and the topological level 
determines the extension up to isomorphism. For any $G$, the levels in the 
image of $\partial$ restrict trivially to $H^3(T)$, and their restriction 
to $H^1(T)\otimes H^2(BT)$ define Weyl-invariant, integral bilinear forms 
on the integer lattice of $T$. From such a form $\beta$, an extension of 
the Lie algebra $L\frg$ is defined by the $2$-cocycle $(\xi,\eta)\mapsto 
\mathrm{Res}_{z=0}\beta(d\xi,\eta)$, and the remaining information about 
the group extension is contained in the torsion part of the level.

\subsection{Adjoint shift.}
A distinguished topological level, the \emph{adjoint shift} $[\sigma]$, 
is the pull-back under the adjoint representation $\mathrm{Ad}:
G \to \SO(\frg)$ of the element $(1,1)$ in $H_{\SO}^3(\SO;\bZ) = \bZ\oplus 
\bZ/2$. The element of order $2$ is pulled back from the integral lift 
$W_3\in H_{\SO}^3(*;\bZ)$ of the third Stiefel-Whitney class, and the 
group is split by restriction to the identity in $\SO$. The free summand 
has a distinguished positive direction, for which the associated form 
$\beta$ is positive definite. When $G$ is simple and simply connected, 
$H_G^3(G;\bZ)\cong H^4(BG\bZ)\cong \bZ$ and $[\sigma] =c$, the dual 
Coxeter number. 

\subsection{Adjoint grading.} 
The adjoint shift also has a component $[\sigma']\in H_G^1
(G;\bZ/2)$, pulled back from the non-trivial class in $H_{\SO}^1
(\SO;\bZ/2) = \bZ/2$. On the loop group side, this gives a homomorphism 
from $LG$ to $\bZ/2$. The presence of this grading-shift means that the 
``usual" Verlinde ring is really the twisted \emph{and} graded, 
${}^{\sigma',\tau}K$-theory; $\xt K$, on the other hand, corresponds to 
a version of the Verlinde group built from the \emph{$\sigma'$-graded 
representations} of the loop group (definitions \ref{A1}, \ref{A2}).

\begin{remark}\label{53} When $\pi_1G$ is free, the most conspicuous 
part of the adjoint shift is the dual Coxeter number. Indeed, $\mathrm 
{Ad}^*W_3=0$; further, if $\tau$ is regular, it turns out that all PERÕs 
of the loop group can be $\sigma'$-graded, so the usual and $\sigma'
$-graded Verlinde groups are (non-canonically) isomorphic. This is 
usually \textit{not} the case when $\pi_1$ has $2$-torsion: see \S\ref
{so3} for $G=\SO(3)$. However, even when $\pi_1$ is free, the $R_G
$-module structures of the $\xt K_G(G)$ are modified by addition of
a grading: the support of the $^\tau\cK(G)$-sheaves is shifted 
(Thm.~\ref{42}.ii; cf.~also Remark \ref{55}). 
\end{remark}

\subsection{Spinors.}The adjoint shift is best understood in terms of 
spinors on $L\frg$. The central extension $\sigma$ of the loop group 
arises by pulling back the (positive energy) Spin representation of
$L\SO(\frg)$; similarly, the grading stems from the \textit{Clifford 
algebra} on $L\frg$. Without unduly labouring this point here, 
we mention that the truly canonical loop group counterpart of $\xt 
K_G^*(G)$ is the $K$-group of graded, $\tau$-projective PERs of the 
crossed product $LG\ltimes \mathrm{Cliff}(L\frg)$. The adjoint and 
dimension shifts arise when relating the latter to PERs of the loop 
group by factoring out the Spin representation. 

\subsection{}By definition, \emph{positive energy} representations of 
$LG$ carry an intertwining \emph{energy} action of the circle group of 
rotations of the loop, with spectrum bounded below. Irreducible PERs have 
\emph{formal characters} in $R_T((q))$, capturing the action of the constant 
loops in $T$, and the energy grading via the variable $q$. We focus on 
the case when $\pi_1G$ is free. On the loop group side, this ensures that 
PERs are determined by their character. The corresponding topological 
fact is injectivity of the restriction $\xt K_G^d(G)\to \xt K_T^d(T)$, 
proved in \S\ref{4}. For the irreducible PER of level $\tau-\sigma$ whose 
zero-energy space has highest weight $\lambda\in\frt^*$, the formal 
character is given by the \emph{Kac formula}
\begin{equation}\label{54}			
\frac{{\sum}_\mu \mathrm{sgn}(\mu)\cdot e^\mu \cdot q^{\|\mu\|^2 -
			\|\lambda +\rho\|^2}}
{\Delta(\frg;q)},
\end{equation}
where $\|x\|^2 = \beta(x,x)/2$ is determined by the level $\tau$, $\mu$ 
ranges over the affine Weyl orbit of $\lambda +\rho$ and $\mathrm{sgn}(\mu)$ 
is the signature of the affine Weyl element taking $\mu$ to $\lambda +\rho$. 
The denominator $\Delta(\frg;q)$ is the formal $(q,T)$ super-character 
of the positive energy spinors on $L\frg/\frt$; this is the \textit
{Kac denominator} multiplied by $e^\rho$. After factoring the (extended) 
affine Weyl group as $W\ltimes\pi_1T$, with $\pi_1T$ mapped to the weight 
lattice by the transpose $H_1(T)\to H^2(BT)$ of $\beta$ and acting there 
by translation, the following is obvious by inspection.

\begin{proposition}\label{56} At $q=1$, the numerator in (\ref{54})  
is the Fourier expansion of a Weyl anti-invariant combination of 
$\delta$-functions on $T$, supported on the regular points of $F$.  
\qed 
\end{proposition}

\begin{remark}\label{55} Two subtleties are concealed in \eqref{54}. 
First, the character is not a function, but a section of the line 
bundle ${}^{\tau-\sigma}\cC h$ of the central extension over $LG$. 
Similarly, the numerator is a section of $\xt\cC h$. The central 
extensions are trivial over $T$, but not canonically so, and the 
weights $\rho,\mu,\lambda$ of $T$ are \textit{projective}, of levels 
$\sigma$, $\tau$ and $\tau-\sigma$. Second, in the presence of an 
$H^1$ component $\sigma'$ of the adjoint shift, formula \eqref{54} 
gives the \textit{super-character} of a graded PER; the usual character 
is obtained by modifying $\mathrm{sgn}(\mu)$ by the $[\sigma']$-grading on 
the lattice $\pi_1T$. 
\end{remark}

\subsection{The $\xt K$-class of a representation.} \label{56b}
There results from Prop.~\ref{56} an obvious isomorphism between 
$\xt K_G^d(G)$ and the complex span of irreducible PERs at level 
$(\tau-\sigma)$: the ``value" of our $K$-class at $f\in F^{reg}$ 
gives the coefficient of the $\delta$-function based there. To read 
off this value, we identify $\xt K_G^d(G)_f^\wedge$ with $\xt H_T^* 
\left(T;\xt\calL(f)\right)$, as in \S\ref{4}, restrict to $\xt H^{top} 
\left(T;\xt\calL(f)\right)$, and integrate over $T$; the answer takes 
values in the fibre of $\xt\cC h$ over $f$. The last integration step 
accounts for the Weyl \textit{anti}-invariance of the answer.

\subsection{Torsion in $\pi_1$.}\label{57}
For more general compact connected $G$, the character formula \eqref{54} 
does not quite determine the representation: much as in the case of 
disconnected compact groups, the non-trivial components of $LG$ carry 
additional information. We did not find the necessary result in the 
literature in the precise form we need,\footnote{\textit{Added in 
revision:} The latest version of Wendt's paper \cite{wend} seems 
to have the needed formulae.} so we only outline the correspondence here, 
hoping to return to it in a future paper. 

A choice of dominant alcove in $\frt$ determines a subgroup $Z'\subset 
LG$ which meets each component of $LG$ whose $\pi_0$-image is torsion 
in a translate of $T$. More precisely, if $N\subset G$ is the normaliser 
of $T$, the subgroup $\Gamma N \subset LG$ of geodesic loops in $N$ is 
an extension of the extended affine Weyl group by $T$, and $Z'$ is the 
subgroup of those torsion components of $\Gamma N$ whose affine co-adjoint 
action preserves the positive alcove. We record, without proof, the 
following fact, reminiscent of the situation of disconnected compact 
groups, where class functions are determined by restriction to a 
quasi-torus.

\begin{proposition}\label{58}
The Kac numerator, at $q=1$, gives a linear isomorphism from the 
complexified space of positive energy representations at level $\tau -
\sigma$ and the space of $\Gamma N$-anti-invariant distributional 
sections of $\xt\cC h$ over $Z'$.
\end{proposition}

\noindent For each component $Z'_w$ of $Z'$, ``anti-invariance" refers 
to its normaliser in $\Gamma N$. Note also that anti-invariant sections 
must be supported on the torsion components of $\Gamma N$: regularity of 
$\tau$ forces the torus to act non-trivially on the fibres of $\xt\cC h$, 
on other components. 

\subsection{} Generalising \eqref{56b}, invariance under the lattice 
part $\pi_1T$ of the affine Weyl group forces a distributional section 
to be a combination of $\delta$-functions supported on isolated $T
$-conjugacy classes in $Z'$. These classes are closely related to the
supports of the ${}^\tau\cK^*(\tilde{G})$ sheaves in the previous
section. Namely, $Z'$ embeds in $G$ via the projection of $\Gamma N$
to $N$ (which kills the lattice in the affine Weyl group). The image 
$Z$ is the largest of all centralisers of dominant regular elements of 
$T$ \eqref{44}. It turns out, thereunder, that the space of $N$-conjugacy 
classes of each component $Z_w$ of $Z$ embeds into the space $T/W$ of 
$G$-conjugacy classes, hence in $Q$. These embeddings furnish a bijection 
between the stalks of the ${}^\tau\cK^*(\tilde{G})$-sheaves for the 
component $Z_w$, in Prop.~\ref{47}, and the invariant $\delta$-sections 
of $\xt\cC h$ supported on $Z'_w$.

\section{The Pontryagin product on $\xt K_G^*(G)$ by localisation}
\label{6}
We now define and compute the Pontryagin product on $\xt K_G^*(G)$,
under the simplifying assumption that $G$ is connected and $\pi_1G$
is free. There is a good \textit{a priori} reason \eqref{63} why the 
answer is a quotient of $R_G$, but we shall make the homomorphism 
explicit by localisation, and recover a well-known description \eqref
{64} of the complex Verlinde algebra. 

In the notation of \S\ref{4}, comparing the sequence $1\to Z_1 \to 
G\to\mathrm{Ad}(G')\to 1$ with the derived sequence $0\to \pi_1G'\to 
\pi_1G\to \pi_1G^{ab}\to 0$ leads to a splitting
\begin{equation}\label{61}		
H_G^3(G) = H_{G'}^3(G')\oplus H^3(G^{ab})\oplus
H^2(BG^{ab})\otimes H^1(G^{ab}).
\end{equation}		
A priori, this holds only rationally, but the absence of torsion 
carries this over to $\bZ$. If $[\tau]$  vanishes in $H^3(T)$, 
then its component in the second summand is null. In this case, $[\tau]$ 
is \textit{equivariantly primitive} for the multiplication map $m:G\times 
G\to G$, meaning that we have an equality, after restriction from 
$H_{G\times G}^3(G\times G;\bZ)$ to the diagonally-equivariant $H_G^3
(G\times G;\bZ)$:
\[
m^*[\tau] = [\tau]\otimes 1+1\otimes[\tau]. 
\]
We denote this restriction by $[\tau,\tau]$. After using the K\"unneth 
and restriction maps 
\[
\xt K_G^*(G)\otimes \xt K_G^*(G)\to {}^{(\tau,\tau)}K_{G\times G}
^*(G\times G)\to {}^{(\tau,\tau)}K_G^*(G\times G),
\]
we define the \textit{Pontryagin} (or \textit{convolution})
\textit{product} on $\xt K_G^*(G)$ as the push-forward along $m$: 
\begin{equation}\label{62}				
m_!:\xt K_G^*(G)\otimes \xt K_G^*(G)\to \xt K_G^*(G).
\end{equation}
The absence of torsion in $H^3$ implies the vanishing of the
Stiefel-Whitney obstruction $W_3$ to $\mathrm{Spin}^c$-orientability, 
so that $m_!$ can be defined. However, choices are involved when
$\pi_1G\ne 0$, and the map (\ref{62}) is only defined up to tensoring with
a 1-dimensional character of $G$. There are two sources for this
ambiguity, and they are resolved in different ways.

The first, and more obvious ambiguity lies in a choice of $\mathrm
{Spin}^c$ structure. This is settled by choosing an Ad-invariant 
$\mathrm{Spin}^c$ structure on $\frg$; but the true explanation of 
this ambiguity and its correct resolution emerge from the use of 
twisted $K$-homology $\xt K_0^G(G)$, in which $m_*$ is the natural map. 
The $\mathrm{Spin}^c$ lifting ambiguity is transferred into the Poincar\'e duality identification $\xt K_0^G(G)\cong\xt K_G^d(G)$ ($d=\dim G$).

The second, and more subtle indeterminacy lies in the twisted $K$-groups
themselves; the class $[\tau]$ only determines $\xt K$ up to tensoring 
by a line bundle. To resolve this ambiguity in (\ref{62}), we need to 
lift the equality $m^*[\tau] = [\tau,\tau]$ from cohomology classes 
to cocycles. When $[\tau]$ is transgressed from an integral class in 
$H^4(BG)$, there is a distinguished such lifting, and a canonical 
multiplication results \cite[I]{fht}. A natural multiplication also 
exists, for any $\tau$, when $\pi_1$ is free, for a different reason. 
As $H_G^3(\{e\})=0$, we can represent $[\tau]$ by a cocycle vanishing 
on $BG\times\{e\}$. This defines an isomorphism $R_G\cong \xt R_G$ and, 
after a choice of $\mathrm{Spin}^c$ structure, a direct image $e_!: R_G
\to \xt K_G^d(G)$. We normalise (\ref{62}) by declaring $m_!(e_!1\otimes 
e_!1) = e_!1$. This is the only choice which makes $\xt K_G^d(G)$ into 
an $R_G$-algebra.

\begin{remark}\label{63} When $\pi_1$ is free, $e_!$ is surjective, 
and $\xt K_G^d(G)$ becomes a quotient ring of $R_G$. There seems to 
be no simple description of the multiplication for general $G$. 
\end{remark}
\begin{theorem}\label{64} When $G$ is connected with free $\pi_1$, 
$\xt K_G^*(G)$ with the Pontryagin product \eqref{62} is isomorphic 
to the ring of Weyl-invariant functions on the regular points of the 
finite subset $F$ of \S\ref{4}, with the point-wise multiplication. 
The isomorphism is the one defined from the Kac numerator (\S\ref{56b}), 
followed by division by the Weyl denominator.

\end{theorem}
\begin{proof} The map \eqref{62} turns $\xt K_G^d(G)$ into an 
$R_G$-algebra. Localising over $Spec(\bC R_G)$, as in (\ref{42}),
we obtain $1$-dimensional fibres over the regular elements of $F/W$. 
We can spell out the ``localised convolution product" as follows. A 
regular point $f\in F$ has centraliser $T$ (Prop.~\ref{43}), and
the localised $K$-theory is isomorphic to $\xt H_T^*(T)$. Dividing 
by the $K$-theoretic Euler class for the inclusion $T\subset G$ 
(the Weyl denominator) converts $G$-convolution to $T$-convolution. 
As noted in (\ref{34}.iii), upon forgetting the $T$-action, the 
generator of $\xt H_T^*(T)$ sees the volume form on $T$; so this 
localised convolution algebra is naturally identified with $\bC$. 
\end{proof} 

\section{Verlinde's formula as a topological index in $\xt K$-theory}
\label{7}
The \textit{Verlinde formula} expresses the dimension of the spaces of
sections of positive holomorphic line bundles over the moduli space $M$
of semi-stable $G_\bC$-bundles over a compact Riemann surface $\Sigma$ of
genus $g>0$. A version of the formula exists for all connected semi-simple 
groups \cite{mein}, but we confine ourselves to simple, simply connected 
ones, in which case $\mathrm{Pic}(M)\cong H^2(M;\bZ)$, line bundles 
$\cO(h)$ being classified by their Chern class $c_1=h$. We shall measure 
$h$ in $H_G^3(G)\cong\bZ$, into which the other two groups embed.\footnote
{The last group can be larger than the others, cf.~\cite{bvl2}; its elements 
classify rank-one \textit{reflexive sheaves} on $M$, and Verlinde's formula 
then expresses their holomorphic Euler characteristic.} Positive line bundles 
have no higher cohomology \cite{kumar}, and the dimension of their space 
of sections is
\begin{equation}\label{71}		
h^0\left(M;\cO(h)\right) = \chi \left(M;\cO(h)\right)=\left| F
\right|^{g-1}\sum\nolimits_{f\in{F^{reg}/W}}{\Delta(f)^{2-2g}},
\end{equation}
where the group $F$ of \S\ref{4} is defined with respect to the shifted 
level $h+c$, and $\Delta$ is the anti-symmetric (spinorial) Weyl denominator. 
We shall replicate the right-hand side of (\ref{71}) in twisted $K$-theory. 
This does not prove the Verlinde formula;\footnote{A simple proof and 
generalisation of the formula, inspired by these ideas, is given in \cite
{tw}.} rather, it interprets it as an infinite-dimensional index theorem, 
in which $\xt K$-theory carries the \emph{topological index}.

\subsection{Reinterpretation of the left side of (\ref{71}).}
Let $\Sigma^\times$ be the complement of a point $z=0$, in a local
coordinate $z$ on $\Sigma$, $G((z))$ the \textit{formal loop group} 
of Laurent series with values in $G_\bC$ and $G[\Sigma^\times]$ the 
subgroup of $G_\bC$-valued algebraic maps on $\Sigma^\times$, and let 
$X:=G((z)) \left/G[\Sigma^\times]\right.$ be a \textit{generalised flag 
variety} for $G((z))$. $X$ is related to $M$ via the quotient stack 
$X/G[[z]]$ by the group of formal regular loops, which is also the 
\textit{moduli stack} of all algebraic $G$-bundles on $\Sigma$. (For
more background on these objects, see \cite{bvl1,falt,sorg,tel}).

\begin{remark}\label{72} If $\Sigma^\circ$ is the complement $\Sigma
-\Delta $ of a small open disk centred at $z=0$, $LG_\bC$ the smooth 
loop group based on $\partial\Delta$ and $Hol(\Sigma^\circ,G_\bC)$ the
subgroup of loops extending holomorphically over $\Sigma^\circ$, then
$X$ is an algebraic model for the homogeneous space $X':=\left.LG_\bC
\right/Hol(\Sigma^\circ,G_\bC)$, which is dense in $X$ and homotopy 
equivalent to it \cite{tel}.
\end{remark}

Algebraic line bundles over $X$ are classified by their Chern class in
$H^2(X;\bZ)\cong H^3(G;\bZ)\cong\bZ$. Theorem 4 of \cite{tel} asserts 
that $H^0(X;\cO(h))$ is a finite sum of duals of irreducible PERs of 
$G((z))$ of level $h$, whereas higher cohomologies vanish. Moreover, 
the multiplicity of the vacuum representation is given by \eqref{71}, 
and more generally, the dual $\mathcal{H}(V)^*$ of the PER with ground 
space $V$ appears with multiplicity
\begin{equation}\label{73}			
m_V=\left| F \right|^{g-1}\sum\nolimits_{f\in {F^{reg}/W}} 
{\Delta(f)^{2-2g}\cdot\chi_V(f)}.
\end{equation}
We can reformulate this, using the inner product\footnote{Its topological 
meaning will be discussed in \cite{fht} I, in connection with the Frobenius 
algebra structure.} on the Verlinde algebra $V(h)$, in which the irreducible 
PERs form an orthonormal basis. On functions on $F^{reg}/W$, this is given 
by
\begin{equation}\label{74}			
\left\langle \varphi \left|\psi\right.\right\rangle =
\left| F \right|^{-1}\sum\nolimits_{f\in F^{reg}/W} 
		{\Delta(f)^2\cdot\bar\varphi(f)\cdot \psi(f)}.
\end{equation}
Formulae (\ref{71}), (\ref{73}) are then captured by 
the following identity in the Verlinde algebra: 
\begin{equation}\label{75}			
\chi\left(X;\cO(h)\right)^* = H^0\left(X;\cO(h)\right)^*
	=\left|F\right|^g\cdot\Delta^{-2g}.
\end{equation}
\vspace{-2ex}
\subsection{Twisted $K$ meaning of the right-hand side.}
Consider the \textit{product of commutators} map $\Pi:G^{2g}\to G$, 
defined by $\Sigma^\circ$: viewing $G^{2g}$ and $G$ as the moduli
spaces of based, flat $G$-bundles over $\Sigma^\circ $, resp.~$\partial
\Sigma^\circ $, $\Pi$ is the restriction of bundles to the boundary.
It is equivariant for $G$-conjugation; dividing by that amounts to
forgetting the base-point, but we shall not do so, and work
equivariantly instead. (Reference to a base-point can be removed by
using the \textit{moduli stacks} of unbased flat bundles, the quotient
stacks $G^{2g}/G$ and $G/G$.) 

The twisting $\tau =h+c\in H_G^3(G)\cong \bZ$, lifted to $H_G^3(G^{2g})$, 
has a distinguished trivialisation, as follows. Any class $[\tau]$ is 
transgressed from a class $[\tau']\in H^4(BG)$, under the classifying 
map $\partial\Sigma^\circ\times G/G\to BG$; in any reasonable model 
for cocycles, $\Pi^*(\tau)$ is the co-boundary of the slant product 
with $\Sigma ^\circ $ of the pull-back of $\tau'$ under $\Sigma^
\circ\times G^{2g}/G\to BG$. So we have a natural isomorphism 
$K_G^*(G^{2g})=\xt K_G^*(G^{2g})$; in particular, a canonical class 
``$\xt 1$" is defined in $\xt K_G^0(G^{2g})$. 

\begin{theorem}\label{76} $\Pi_!(\xt 1)\in \xt K_G^d(G)$ is the 
right-hand side of (\ref{75}). In other words, the multiplicities of 
$\Pi_!(\xt 1)$ in the basis of PERs are the Verlinde numbers (\ref{73}). 
\end{theorem}

\begin{proof} Localising to a regular $f\in F/W$, $\Pi_!(\xt 1)$ 
agrees with the push-forward of $\xt 1/\Delta(f)^{2g}$ along $T^{2g}$, 
because $\Delta(f)^{2g}$ is the relative $K$-theory Euler class of 
the embedding $T^{2g}\subset G^{2g}$. $T^{2g}$ maps to $e\in G$, 
so $\Pi_!$ factors as the push-forward to $K_T^0(*)$, followed 
by $e_!$ of \S\ref{6}. It might seem, at first glance, that the 
$K$-theoretic integral of $\xt 1$ over $T^{2g}$ is nought, but that 
is not so. In factoring $\Pi_!$, we have trivialised $\Pi^*\tau$ 
on $T^{2g}$ by doing so first in $H_G^3(e)$. This differs from the 
trivialisation by $\Sigma^\circ$-transgression, the difference being 
the $2$-dimensional transgression over $\Sigma$ of the pull-back of 
$\tau'$ under the classifying map $\Sigma\times T^{2g}\to BT$. The
difference line bundle over $T^{2g}$ is the restriction of $\cO(h+c)$,
and its integral over $T^{2g}$ is $\left|F\right|^g$ \cite{mein}. 
All in all, the $T$-restriction of our push-forward is $\left|F\right|
^g\cdot\Delta(f)^{-2g}\cdot e_!(1)$; but this is the value at $f$ of 
the right side of (\ref{75}).  
\end{proof}

\subsection{Interpretation as an index theorem.}
We have the following set-up in mind. Let $\pi:S\to B$ be a proper
submersion of manifolds of relative dimension $d$, $\tau$ a twisting 
of $K$-theory over $B$ whose lifting to $S$ is trivial. Expressing 
$\pi^*\tau $ as a co-boundary $\delta \omega$ on $S$ fixes an 
isomorphism ${}^{\pi^*\tau}\!K^*(S)\cong K^*(S)$. As a result, a 
class $\xt 1\in {}^{\pi^*\tau}\!K^*(S)$ is unambiguously defined. 
If $\pi$ is $K$-oriented, we obtain an index class $\pi_!\xt 1
\in\xt K^d(B)$. 

On any fibre $S_b$, $\pi^*\tau$, being lifted from a point, is null 
as a cocycle, which again allows one to identify ${}^{\pi^*\tau}\!
K^*(S_b)$ with $K^*(S_b)$. However, the new identification differs from 
the old one $\pi^*\tau =\delta\omega $; so the restriction of $\xt 1$ 
corresponds to a line bundle $\cO(\omega)$ over $S_b$, with $c_1=\omega$. 
The fibre of the relative index bundle $\pi_!\xt 1$ at $b$ is then the 
index of $\cO(\omega)$ over $S_b$.

In our situation, the index (\ref{75}) ought to be captured by the 
map from the manifold $X'$ (\ref{72}) to a point, in $LG$-equivariant 
$K$-theory. Take for $B$ the \textit{classifying stack} $BLG$, and for 
$S$ the \textit{quotient stack} $X'/LG$. As a real manifold, $X'$ is 
the moduli space of flat $G$-connections on $\Sigma ^\circ $, trivialised 
on the boundary; so $X'/LG$ is equivalent to the quotient stack $G^{2g}/G$. 
$BLG$ has the homotopy type of $G/G$, and in these identifications, the 
projection to a point $X'\to *$ becomes our map $\Pi:G^{2g}/G\to G/G$.  

Twistings are required, since the $LG$-action is projective on $H^0$ and
on the line bundle. In \S\ref{5}, we identified ${}^{h+c}K_G^d(G)$ with 
the Verlinde algebra $V(h)$, which is where our analytic index lives. 
This leads us to the push-forward $\Pi_!:{}^{h+c}K_G^0(G^{2g})\to
{}^{h+c}K_G^0(G)$. We can reconcile this shifted twisting with the
level $h$ in \eqref{75} by reinterpreting the left-hand side there as the
\textit{Dirac index} of $(h+c)$; the $c$-shift is the projective cocycle
of the $LG$-action on spinors on $X$. (This is one of the ways in which
$(-2c)$ behaves like the canonical bundle of $X$.) Thus, \eqref{75} and 
\eqref{76} express the equality of the algebraic and topological indexes.

\subsection{Summary.} 
We summarise and clarify our discussion by introducing the space $\cA$ 
of $G$-connections on $\partial\Sigma^\circ $. The stack 
$\cA/LG$ is equivalent to $G/G$; so ${}^{h+c}K_G^d(G)={}^{h+c}K_{LG}^d
(\cA)$; moreover, the boundary restriction $X'/LG \to \cA/LG$ is exactly 
$\Pi$. We are then saying that the $LG$-equivariant analytic index over 
$X'$, rigorously defined in the algebraic model, is computed topologically, 
by factoring the push-forward to a point into the (rigorously defined) 
Gysin map $\Pi_!:{}^{h+c}K_{LG}^0(X')\to {}^{h+c}K_{LG}^0(\cA)$, along
the proper map $\Pi$, and the $LG$-equivariant push-forward along 
$\cA\to *$; the latter, we interpret as the isomorphism 
${}^{h+c}K_G^d(G)\cong V(h)$. 

\section{A worked example: $G = \SO(3)$}\label{so3}
We now illustrate the effect of gradings on the twisted $K_G(G)$ for the 
group $G=\SO(3)$. This is simple enough to settle integrally, by an 
independent Mayer-Vietoris calculation. There is no true difficulty 
in extending the Chern character calculation to $\SO(n)$, as in Example 
\ref{46}, but the integral treatment requires other methods \cite{fht}.

\subsection{Twistings.} \label{B1} $H^1_G(G;\bZ/2) = \bZ/2$, representing 
the double cover $\SU(2)$, while $H^3_G(G;\bZ) = \bZ/2\oplus\bZ$, split by 
the restriction to $H^3_G(\{e\}) = \bZ/2$. We label a twisting class, 
inclusive of the grading, by a triple in $\{\pm\}\times\{\pm\}\times\bZ$, 
ordered as introduced: $H^1,H^3_{tors}, H^3_{free}$. For instance, the \textit{adjoint shift} $[\sigma'+\sigma]$ of \S\ref{5} is $(-,-,1)$.

The even elements in $H^3_G(G;\bZ)$ with no torsion component are 
transgressed from $H^4(BG;\bZ)$, and hence are primitive classes (\S\ref
{6}). On the other hand, the \textit{transgressed $K$-theory twistings} 
are the types $(+,+,\mbox{even})$ and $(-,+,\mbox{odd})$. This stems 
from the fact that the classifying space for the group of twistings 
is a \textit{non-trivial} fibration\footnote{This leads to a non-standard 
group structure on twisting classes, which is the extension of $H^1(\bZ/2)$ 
by $H^3(\bZ)$ using the Bockstein of the cup-product; cf.~\cite
{atseg2,fht}.} $K(\bZ;4)\to E \to K(\bZ/2;2)$, with $k$-invariant 
the integral Bockstein of $Sq^2$; we have $[BG;E]\cong \bZ$, with odd 
classes mapping to the non-trivial element of $H^2(BG;\bZ)$. With 
reference to our discussion in \S\ref{6}, and correcting for the 
Ad-equivariant $W_3(G)$, it follows that $K_G(G)$ carries a natural 
Pontryagin product for the twisting types $(+,-,\mbox{even})$ and $(-,-,
\mbox{odd})$.

\subsection{Spin structures.} One fine point that we will not fully explore 
here concerns the r\^ole of Spin structures. In the presence of an $H^2$ 
component, integrating an $E$-valued class requires a choice of Spin 
structure. The class $(-,+, \mbox{odd})$ is transgressed from an (odd) 
class in $[BG;E]$ using the \textit{bounding} (non-trivial) Spin structure 
on the circle, which is the one leading to a Pontryagin product, by extension 
to the three-holed sphere; transgressing over the non-bounding structure 
leads to $(-,-, \mbox{odd})$. The $K$-group with $W_3$-shifted twisting 
$(-,+, \mbox{odd})$ is a module over the other $K$-ring, and indeed the 
two $K$-groups assemble to a \textit{field theory for Spin surfaces}.       

\begin{remark}\label{B2}\begin{trivlist}\itemsep0ex
\item (i) Shifting by $[\sigma]$, it follows that the positive energy 
representations of $L\SO(3)$ will carry a fusion ring structure 
for levels of type $(+,+,\mbox{even})$ and $(-,+,\mbox{odd})$. The 
parity restriction on the levels was known (G.~Segal, unpublished). 
However, the possibility of defining a fusion product on \textit{graded} 
representations at \textit{odd} level  seems new.
\item (ii) In addition to the two torsion components, we will find that 
the parity of the free component of the twisting affects the form of the 
answer, leading to eight separate cases. Note that elements of $H^3\left
(\SO(3);\bZ\right)$ lift to the even classes in $H^3\left(\SU(2);
\bZ\right)$, so parity distinguishes the residues $0$ and $2\pmod{4}$ 
of the level in $\SU(2)$. 
\end{trivlist}
\end{remark}

\subsection{The $\cK$-sheaves.} We have $G = \SU(2)/\!\pm\!\mathrm{I}_2$, 
and, with reference to Example \ref{25}, the maximal torus $T_\bC=\SO(2;
\bC)$ is parametrised by $\mu=\lambda^2\in\bC^\times$; the Weyl group 
interchanges $\mu$ with $\mu^{-1}$. The centraliser $Z(\mu)$ is $T$, 
except that $Z(1)=\SO(3)$ and $Z(-1)=\Or(2)$. With $\tau=(\vep_1,
\vep_2,k)$, $k>0$, the twisted cohomologies vanish for $\mu=1$, as in 
(\ref{34}). The holonomy of $\xt\calL(\mu)$ around $T$ is $\vep_1\mu^k$:
to see this, use Remark \ref{B2}.ii, noting that $T$ is half the size of 
its counterpart in (\ref{25}). We thus find $1$-dimensional stalks of 
$^\tau\cK^1_G$ at every $k$th root of $\vep_1$ in the upper-half plane, 
and no contribution at $\mu=1$. Further investigation is needed for 
$\mu=-1$. 

We need the action of $Z(-1)=\Or(2)$ on $\xt H^*_T\left(\Or(2);\xt\calL
(-1)\right)$. Now, $\bT$-central extensions of $\Or(2)$ are classified 
by $\bZ/2$, and the twisting $\tau$, restricted to $H^3(B\Or(2);\bZ)$, 
can be traced by Mayer-Vietoris to define the extension class $(-1)^k
\vep_2$ and the grading $\vep_1$. In particular, $\Or(2)$ acts on the 
fibre of $\xt\calL(-1)$ over $(-1)$ via the character $(-1)^k\vep_1\vep_2$ 
(cf.~Construction \ref{212}), hence by the opposite character on $H^1
(T;\xt\calL)$. Further, the restriction to the $\Or(2)$-stabiliser $\langle
w,(-1)\rangle\cong\bZ/2\times\bZ/2$ of a point $w$ in the other component 
of $Z(-1)$ is $(-1)^k\vep_2$; moreover, $\mu=-1$ will have parity $\vep_1$ 
for the odd component of $Z(-1)$ (cf.~\S\ref{A14}). The stabiliser will 
therefore act by $w\mapsto(-1)^k\vep_2$, $(-1)\mapsto 1$ on the fibre of 
$\xt\calL(-1)$ and on the cohomology $H^0_T\left(Z(-1)_w;\xt\calL\right)$. 

Therefore, at $\mu=-1$, we will get a $1$-dimensional contribution to 
$^\tau\cK^1_G$ from $T$ when $\tau$ has type $(+,-, \mathrm{even})$ and 
$(-,-,\mathrm{odd})$. We get, by Thm.~\ref{A14}, a $1$-dimensional 
contribution to $\cK^0$ from the odd component of $Z(-1)$ for the twistings 
$(+,+,\mathrm{even})$ and $(+,-,\mathrm{odd})$, and a $1$-dimensional 
contribution to $\cK^1$ for $(-,+,\mathrm{even})$ and $(-,-,
\mathrm{odd})$.

\subsection{} The representation ring $R_{\SU(2)}$ splits into $^+R:=
R_{\SO(3)}$ and its $H^3$-twisted version $^-R$. The quotient $R(k)$ of $R$ 
by the principal ideal on the irreducible $2k$-dimensional representation, 
spanned over $\bZ$ by the images $[1],\ldots,[2k-1]$ of the irreducible 
representations of those dimensions, splits similarly as $^+R(k)$ and 
$^-R(k)$. Denote by ${}^{\vep_1,\vep_2} R(k)$ ($\vep_i\in\{\pm1\}$) the 
quotient of $^{\vep_2}R(k)$ under the identifications $[p]\sim\vep_1[2k-p]$ 
(and also $[k]=0$ if $\vep_1=-$). The quotients $^{\pm,+}R(k)$ turn out to 
be $^+R$-algebras, while $^{\pm,-}R(k)$ are $^+R$-modules. 

\begin{proposition} \label{B4}
The abelian group $K^*_G(G;\bZ)$ is given by the following table; in all 
cases, the $R(k)$ summand is the direct image of $^{\vep_2} R_{\mathrm 
{SO}(3)}$ under the inclusion $\{e\}\hookrightarrow G$.
\vskip1ex
\begin{tabular}{lcll}
$\stackrel{\mathrm{Twisting}}{(k\,\mathrm{even})}$ && $K^0$ & $\quad K^1$ \\
\hline\hline
$(+,+,k)$  &\vline& $\bZ$ & $^{--}R(k)$\\ 
$(-,+,k)$  &\vline& $0$   & $^{+-}R(k)\oplus\bZ$\\ 
$(+,-,k)^*$  &\vline& $0$   & $^{-+}R(k)$ \\ 
$(-,-,k)$  &\vline& $0$   & $^{++}R(k)$ \\
\hline
\end{tabular}
\hspace{2.5em}
\begin{tabular}{lcll}
$\stackrel{\mathrm{Twisting}}{(k\,\mathrm{odd})}$ && $K^0$ & $\quad K^1$\\
\hline\hline
$(+,+,k)$ &\vline& $0$ & $^{--}R(k)$ \\
$(-,+,k)$ &\vline& $0$ & $^{+-}R(k)$\\
$(+,-,k)$ &\vline& $\bZ$ & $^{-+}R(k)$\\
$(-,-,k)^*$ &\vline& $0$ & $^{++}R(k)\oplus\bZ$\\
\hline
\end{tabular}
\end{proposition}
\noindent The starred twistings lead to a natural Pontryagin product on 
$K$-theory.

\begin{proof} In order to identify the induction map $R\to K_G(G)$, we 
use a Mayer-Vietoris calculation in \textit{$K$-homology}; Poincar\'e 
duality along $G$ gives an isomorphism 
\[
{}^{\vep_1,\vep_2,k}K^*_G(G) \cong {}^{\vep_1,-\vep_2,k}K_{3-*}^G(G).
\]
Cover $\SO(3)$ by the open sets $U_0, U_\pi$ of rotations by angles different from $\pi$ and $0$, respectively. They deform Ad-invariantly 
onto the conjugacy classes $\{1\}$ and $\bR\bP^2$ (of rotations by angle 
$\pi$), respectively, and their intersection is $G$-equivalent to the 
sphere. The equivariant $K$-homologies in the MV sequence are twisted 
versions of the representation rings of the stabilisers $\SO(3)$, 
$\Or(2)$ and $\SO(2)$. After tracing the twistings, the 
vanishing of all twisted equivariant $K^1$ for the first two groups turns 
the six-term Mayer-Vietoris sequence into the following, where the sign 
$\pm=(-1)^k$:
\begin{equation}\label{B5}
0 \to{}^{-\vep_1,\pm\vep_2}R^1_{\Or(2)} \to \xt K_1^G(G)
\to {}^{\vep_2}R_{\SO(2)} \xrightarrow{\mathrm{Ind}} 
{}^{-\vep_2}R_{\SO(3)}\oplus{}^{-\vep_1,\pm\vep_2}
R_{\Or(2)} \to \xt K_0^G(G) \to 0.
\end{equation}
The first component of Ind is Dirac induction from $\SO(2)$ to 
$\SO(3)$, sending $u^{\pm p}$ to the irreducible virtual 
representation $\pm\bC^{2p}$ ($p$ is a positive integer or half-integer,
according to $\vep_2$). The second component is induction, \textit{preceded 
by tensoring by $u^{-k/2}$}, which is the effect of the twisting.
When $k\neq0$, the total Ind is injective, so (\ref{B5}) breaks into
\begin{align}\label{B7}
\xt K_1^G(G)&\cong {}^{-\vep_1,\pm\vep_2}R^1_{\Or(2)},\nonumber\\
\xt K_0^G(G)&\cong\left.\left({}^{-\vep_2}R_{\SO(3)}\oplus{}^
{-\vep_1,\pm\vep_2}R_{\Or(2)}\right)\right/
{}^{\vep_2}R_{\SO(2)}. 
\end{align}
The second component of Ind is surjective, except when $-\vep_1 =1$ and
$\vep_2=(-1)^k$, in which case the target $R_{\Or(2)}$ has a extra $\bZ$, 
from the sign representation. So, except in that case, $K^1_G(G)$ is a 
quotient of $^+R$; we leave it to the reader to identify it with the 
appropriate version of $R(k)$, as in Prop.~\ref{B4}. Similarly, the $R^1$ 
group is nil, except when $\vep_1 = 1$ and $\vep_2=(-1)^k$, in which case 
it is $\bZ$ (spanned by the graded, $2$-dimensional representation of 
$\bZ/2$). \end{proof}

\begin{remark}
Because induction from $\SO(2)$ to $\SO(3)$ is always surjective, $K^G_*(G)$ 
is always isomorphic, qua $R_{\SO(3)}$-module, to a quotient of the 
appropriately twisted $R^*_{\Or(2)}$. 
\end{remark}

\subsection{The ring structure.} For transgressed twistings (\S\ref{B1}), 
the Pontryagin product gives an $R_G$-algebra structure on $\xt K_G(G)$.
When $\tau=(+,-,\mbox{even})$, the latter is a quotient of $R_G$ and the 
ring structure is thus completely determined; but for $\tau=(-,-,\mbox{odd})$, 
we need more information. In general, this requires us to compute the 
convolution on $\xt K_Z(Z)$, for the disconnected centralisers; but a 
shortcut is possible for $\SO(3)$. 

Let $[k]_\pm$ be the elements of $\xt K_G(G)$ induced from the $\pm$-sign
representation of $\Or(2)$. From \eqref{B7}, we see that $[k]=[k]_++[k]_-$, 
and that $\{[1],\ldots,[k-2],[k]_+,[k]_-\}$ is an integral basis of the 
Verlinde ring. Multiplication by every $[p]$ is known from the $R_G$-module 
structure; e.g.~, for $p\leq q\leq k$, $[p]\cdot[q] = [q-p+1]+\ldots+[q+p-1]$,
identifying $[n]=[2k-n]$ for $n>k$. Moreover, 
\[
[p]\cdot[k]_\pm = [k-p+1]+\ldots+[k-2]+[k]_{\pm i^{(p-1)}};
\]
note, indeed, that $[k]_\pm=\mbox{Ind}(\bC_\pm)$, induced from the sign
representation of $\Or(2)$, so our product is $\mbox{Ind}(\bC^p\otimes
\bC_\pm)$, and the final $[k]$-sign matches the action of $\Or(2)$ on the 
zero-weight space therein. The missing information involves separating the 
two terms in the sum 
\begin{equation}\label{B10}
[k]_+^2+[k]_+[k]_- = [k][k]_+ =[1]+\ldots+[k-2]+[k]_{i^{(k-1)}}.
\end{equation}
We now observe, in the localised picture of the complex Verlinde algebra, 
that $[k]_+=[k]_-$ \emph{away from} $\mu=-1$. We must then split up the 
right-hand side of \eqref{B10} into two representations whose characters
agree at all $\mu\neq -1$ with $\mu^k=-1$. One such splitting is into residue 
classes $\mod{4}$, and there is no other, by reason of degree (see the 
remark that follows). We then have
\begin{align*}
[k]_+^2 &= \sum_{0<p<k/4}[k-4p] + [k]_{i^{(k-1)}}\\
[k]_+\cdot[k]_- &= \sum_{0<p+1/2<k/4}[k-4p-2]\\
[k]_-^2 &= \sum_{0<p<k/4}[k-4p]+[k]_{-i^{(k-1)}}
\end{align*}
The formulae for $[k]_+^2$ and $[k]_+[k]_-$ could not be switched,
as $[k]_{i^{(k-1)}}$ does \emph{not} appear in $[k][k]_-$.

\begin{remark} We implicitly used the positivity of the ``fusion product" 
coefficients in our argument; without that, other splittings of \eqref{B10} 
are possible, and could only be ruled out by computing the convolution on
$\Or(2)$. We do not know a \emph{topological} argument for positivity; 
rather, this results form a cohomology vanishing theorem on the loop 
group side.
\end{remark}

\section*{Acknowledgments} 
We are indebted to Graeme Segal for helpful conversations 
on foundational aspects of twisted $K$-theory. We also thank the referee 
for flagging several small errors in our original submission and 
suggesting improvements to the exposition. D.S.F and C.T. were 
partially supported by NSF grant DMS-0072675; M.J.H. was partially 
supported by NSF grant DMS-9803428. C.T. also thanks MSRI, 
where much of this work was completed, for its hospitable and 
stimulating environment.

\vskip1cm
{\small
\noindent
DSF: Department of Mathematics, UT Austin, Austin, TX 78712, USA \\
\texttt{dafr@math.utexas.edu}\\
MJH: Department of Mathematics, Harvard University, Cambridge, 
MA 02138, USA \\
\texttt{mjh@math.harvard.edu} \\
CT: DPMMS, CMS, Wilberforce Road, Cambridge CB3 0WB, UK \\
\texttt{teleman@dpmms.cam.ac.uk}}

\end{document}